\documentclass{siamltex}

\usepackage{amssymb}
\usepackage{amsmath}
\usepackage{color}
\usepackage[dvips]{graphicx}
\usepackage{epsfig,subfigure}
\usepackage[makeroom]{cancel}

\DeclareMathOperator*{\argmin}{arg\,min}				   
\renewcommand{\t} {^{\top}}								   
\newcommand{\norm} [2][]{\left\|#2\right\|_{#1}}		   
\newcommand{\trace}{\mathsf{trace}}

\newcommand{\bfGamma}{{
{\Gamma}}}
\newcommand{\bfDelta}{{\boldsymbol{\Delta}}}

\newcommand{\bfLambda}{{\boldsymbol{\Lambda}}}

\newcommand{\bfPhi}{{\boldsymbol{\Phi}}}
\newcommand{\bfPsi}{{\boldsymbol{\Psi}}}

\newcommand{\bfgamma}{{\boldsymbol{\gamma}}}

\newcommand{\bfeta}{{\boldsymbol{\eta}}}

\newcommand{\bflambda}{{\boldsymbol{\lambda}}}

\newcommand{\bfxi}{{\boldsymbol{\xi}}}

\newcommand{\bfphi}{{\boldsymbol{\phi}}}

\newcommand{\bfpsi}{{\boldsymbol{\psi}}}

\newcommand{\bfA}{{\bf A}}

\newcommand{\bfC}{{\bf C}}
\newcommand{\bfD}{{\bf D}}

\newcommand{\bfF}{{\bf F}}

\newcommand{\bfH}{{\bf H}}
\newcommand{\bfI}{{\bf I}}
\newcommand{\bfJ}{{\bf J}}

\newcommand{\bfL}{{\bf L}}

\newcommand{\bfP}{{\bf P}}
\newcommand{\bfQ}{{\bf Q}}
\newcommand{\bfR}{{\bf R}}
\newcommand{\bfS}{{\bf S}}
\newcommand{\bfT}{{\bf T}}

\newcommand{\bfW}{{\bf W}}

\newcommand{\bfZ}{{\bf Z}}

\newcommand{\bfb}{{\bf b}}

\newcommand{\bfe}{{\bf e}}

\newcommand{\bfg}{{\bf g}}

\newcommand{\bfn}{{\bf n}}
\newcommand{\bfp}{{\bf p}}
\newcommand{\bfq}{{\bf q}}

\newcommand{\bfx}{{\bf x}}
\newcommand{\bfy}{{\bf y}}
\newcommand{\bfz}{{\bf z}}
\newcommand{\bfzero}{{\bf0}}

\newcommand{\mxL}{{\bf L}}

\newcommand{\bbR}{\mathbb{R}}

\newcommand{\GCV}{\mathrm{GCV}}
\newcommand{\true}{\mathrm{true}}

\newcommand{\bvec}{\left[ \begin{array}{c} }
\newcommand{\evec}{\end{array} \right]}

\title{Optimal Regularization Parameters for General-Form Tikhonov Regularization}

\author{ Julianne Chung\thanks{Department of Mathematics, Virginia Tech, Blacksburg, VA 24061}
\and Malena I. Espa\~ nol\thanks{Department of Mathematics, The University of Akron, Akron, OH 44325}
\and Tuan Nguyen$^*$}

\begin{document}
\maketitle

\begin{abstract} In this work we consider the problem of finding optimal regularization parameters for general-form Tikhonov regularization using training data.
We formulate the general-form Tikhonov solution as a spectral filtered solution using the generalized singular value decomposition of the matrix of the forward model and a given regularization matrix. Then, we find the optimal regularization parameter by minimizing the average of the errors between the filtered solutions and the true data. We extend the approach to the multi-parameter Tikhonov problem for the case where all the matrices involved are simultaneously diagonalizable. For problems where this is not the case, we describe an approach to compute optimal or near-optimal regularization parameters by using operator approximations for the original problem. Several tests are performed for 1D and 2D examples using different norms on the errors, showing the effectiveness of this approach.
\end{abstract}

\begin{keywords}{spectral filtering, regularization, multi-parameter Tikhonov, optimal filters, learning approach}\end{keywords}

\section{Introduction} 
\label{sec:introduction}

We consider a linear ill-posed problem that can be modeled as
\begin{equation} \label{OrigProb}
	\bfA \bfx_\true +\bfn= \bfb,
\end{equation}
where $\bfx_\true \in \mathbb{R}^{n\times 1}$ is the desired solution, $\bfA\in \mathbb{R}^{m\times n}, m \geq n,$ models the forward process, $\bfn \in \mathbb{R}^{m\times1}$ is additive random noise, and $\bfb \in \mathbb{R}^{m\times1}$ is observed data. The goal is to obtain an approximate solution of $\bfx_\true$ given $\bfA$ and $\bfb$.  Prior information regarding the probability distribution of $\bfn$ may be incorporated, if known. The matrix $\bfA$ is ill-conditioned and its singular values decay to zero without significant spectral gap. Due to the ill-posed nature of the problem, the exact (inverse) solution of \eqref{OrigProb} will be contaminated by noise and may not be a good approximation of $\bfx_\true$. Thus, we seek an approximation to the solution by solving a nearby problem that is well-posed, a process called \emph{regularization}. A well-known regularization method is Tikhonov regularization~\cite{HansenBook1998},
\begin{equation}\label{GeneralTik}
  \bfx_\lambda =
  \argmin_\bfx \left\{ \|\bfA \bfx - \bfb \|_2^2 + \lambda^2 \|\bfL \bfx\|_2^2 \right\},
\end{equation}
where $\lambda$ is a regularization parameter and $\bfL$ is a regularization matrix. The matrix $\bfL$ must have the same number of columns as $\bfA$, and we assume that the null spaces of $\bfA$ and $\bfL$ intersect trivially such that the matrix $[\bfA\t\, \bfL\t ]$ has full row rank and the
solution $\bfx_\lambda$ is unique. If $\lambda$ and $\bfL$ are chosen appropriately, the solution to~\eqref{GeneralTik} should approximate the desired solution $\bfx_\true$. Typical choices of $\bfL$ include the identity matrix and the discrete first or second derivative operators. If $\bfL = \bfI$, we say that \eqref{GeneralTik} is in \emph{standard form}. Otherwise, we say it is in \emph{general form}. For the latter case, we can transform the problem to standard form. That is,  we can use the substitution
$\bfx = \bfL^{-1}\bfy$ if $\bfL$ is invertible, or otherwise, we can define $\bfx = \bfL^\dagger_\bfA \bfy$, where $\bfL^\dagger_\bfA$ is the $\bfA$-weighted generalized inverse of $\bfL$ defined by $\bfL^{\dagger}_\bfA = (\bfI - ( \bfA (\bfI - \bfL^\dagger \bfL))^\dagger \bfA)\bfL^\dagger$; see~\cite[Section~2.3]{HansenBook1998} for further details.

\begin{figure}
	\begin{center}
		\includegraphics[scale=0.65]{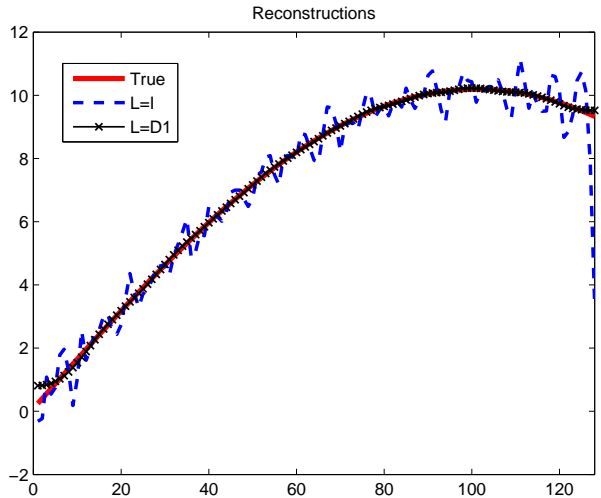}
	\end{center}
	\caption{Motivating example to illustrate the effect of using different regularization matrices $\bfL$ on the reconstructed solution.  Here, the true solution is provided, along with reconstructions corresponding to $\bfL=\bfI$ (standard Tikhonov) and $\bfL=\bfD1$, where $\bfD1$ represents a discrete approximation of the first derivative.}
	\label{fig:differentL}
	\end{figure}

Next we provide a motivating example to illustrate the potential benefits of considering the general-form Tikhonov regularization. Using the deriv2 test problem from \cite{hansen}, we provide in Figure~\ref{fig:differentL} reconstructions for various choices of $\bfL,$ where in each case, the regularization parameter $\lambda$ was selected to produce the smallest mean squared error.
	This example illustrates that $\bfL = \bfI$ is not suited for all problems and that better reconstructions may be obtained with $\bfL \neq \bfI$.  We remark that this example is used for illustrative purposes only, so specific details are not included here. Complete numerical investigations will be presented in Section~\ref{sec:numerics}.
	
The difficulty in selecting an appropriate regularization matrix has recently led researchers to investigate an extension of the general-form Tikhonov regularization that contains several regularization terms instead of only one. This extension is known as the multi-parameter Tikhonov problem,
 \begin{equation}\label{MultiTik}
\bfx_\bflambda = \argmin_\bfx \left\{ \|\bfA \bfx - \bfb \|_2^2 + \sum_{j=1}^J\lambda_j^2 \|\bfL_j \bfx\|_2^2 \right\},
 \end{equation}
where multiple regularization parameters $\bflambda = [\lambda_1,\dots,\lambda_J]$
and multiple regularization matrices $\bfL_j \in \mathbb{R}^{p_j\times n}$ are incorporated.  Several studies on the multi-parameter Tikhonov problem include \cite{BazanBorgesFrancisco2012,BelgeKilmerMiller,BrezinskiMulti2003,GazzolaNovati2013,LuPereverzev2011,Wang2012}.

A key ingredient for both the one-parameter~\eqref{GeneralTik} and multi-parameter~\eqref{MultiTik} Tikhonov problems is the choice of the regularization parameter(s). Standard techniques to find the regularization parameter for Tikhonov regularization have been studied extensively in the literature and include the discrepancy principle, the L-curve criterion, and the generalized cross-validation (GCV) method~\cite{HansenBook1998}. Some of these techniques have been extended to the multi-parameter Tikhonov problem. For example, a higher dimension {L}-curve is considered in \cite{BelgeKilmerMiller}, the discrepancy principle is discussed in \cite{GazzolaNovati2013,LuPereverzev2011,Wang2012}, and an extension of GCV is described in \cite{BrezinskiMulti2003}. In \cite{BazanBorgesFrancisco2012}, the Regi\'nska's parameter choice rule is generalized to the multi-parameter Tikhonov case.

In this paper, we propose a learning approach for computing \emph{optimal} regularization parameters for both Tikhonov problems~\eqref{GeneralTik} and~\eqref{MultiTik}. Previous work on learning approaches in the context of regularization methods for solving inverse problems can be found in~\cite{chung2011designing,ChungChung2013,de2013image,haber2003learning,horesh2009sensitivity,huang2012optimal,kunisch2013bilevel,peyre2011learning}. However, none of these works specifically address the general-form Tikhonov and multi-parameter Tikhonov problems.

We follow the work of Chung et.~al.~\cite{chung2011designing} where training data is used to find optimal regularization parameters for a variety of spectral regularization methods, including standard form Tikhonov regularization.  We first extend these approaches to the general-form Tikhonov problem.
For general-form Tikhonov, the generalized singular value decomposition (GSVD)~\cite{HansenBook2010} of  $\{\bfA, \bfL\}$ can be used to write the solution to \eqref{GeneralTik} as a spectral filtered solution.  We take advantage of this fact to develop efficient methods for computing an optimal regularization parameter, and compare its performance to that of filtered SVD solutions on several numerical examples.

Another main contribution of our work is to develop an efficient learning approach for the multi-parameter Tikhonov problem.
The extension to the multi-parameter case is more complex, not only in terms of finding multiple
$\lambda$'s, but also because in general, there is no natural extension of the GSVD for
 $\bfA$ and multiple regularization matrices $\bfL_j$. For problems where matrices $\bfA$ and $\bfL_j$ are \emph{simultaneously diagonalizable} (SD), such as those arising in image processing, we propose an efficient learning approach for computing optimal regularization parameters.  However, for problems that do not exhibit this nice property, we follow
a framework for regularization introduced in~\cite{ChungKilmerOleary}, where matrix approximations are used to compute regularization parameters, but the original problem is solved.
It is worth mentioning that in~\cite{de2013image} and \cite{kunisch2013bilevel}, multi-parameter learning approaches for denoising problems ($\bfA=\bfI$) are proposed. In~\cite{kunisch2013bilevel} a parameter learning approach for multiple $\ell_p$-norm regularization terms is presented. The learning problem is formulated as a bilevel optimization problem and solved using semismooth Newton methods. Although the same techniques could be extended to our problem, we believe that such an approach applied to the problem where $\bfA\neq\bfI$ would be more computationally expensive than the one we propose here.

The paper is organized as follows.  Background on the general-form Tikhonov regularization problem is provided in Section \ref{sec:background}.  We provide the solution as a spectral filtered solution and describe some standard approaches for selecting regularization parameters. In Section \ref{sec:optfilt}, we extend the optimal filter framework to both the general-form Tikhonov regularization problem and the multi-parameter Tikhonov regularization problem and describe an empirical Bayes risk framework for computing optimal regularization parameters. In particular, numerical methods for computing optimal regularization parameters for the multi-parameter problem are discussed. Numerical results are presented in Section \ref{sec:numerics}, and conclusions and future work can be found in Section \ref{sec:conclusions}.

\section{Background} 
\label{sec:background}
A closed form solution to the general-form Tikhonov regularization problem~\eqref{GeneralTik} can be obtained using the GSVD of the matrix pair $\{\bfA,\bfL\}$.
Given matrices $\bfA \in \bbR^{m \times n}, m \geq n$, and $\bfL \in \bbR^{p \times n}$, the GSVD can be obtained by first considering the {\it reduced} QR-factorization of the stacked matrix,
\[\left[\begin{array}{c} \bfA \\ \bfL \end{array}\right] = \left[\begin{array}{c} \bfQ_{\bfA} \\ \bfQ_{\bfL} \end{array}\right] \bfR . \] We have
$\bfA=\bfQ_{\bfA}\bfR$ and $\bfL=\bfQ_{\bfL}\bfR$. Let $\bfQ_{\bfA}=\bfP \bfC \bfW^{T}$ and
$\bfQ_{\bfL}=\bar{\bfP}\bfS \bfW^{T}$ be the CS decomposition \cite{GolubVanLoanBook} of
\{$\bfQ_{\bfA},\bfQ_{\bfL}$\}, where $\bfP\in\mathbb{R}^{m\times m}$,
$\bar{\bfP}\in\mathbb{R}^{p\times p}$, and $\bfW\in\mathbb{R}^{n\times n}$
are orthogonal matrices; $\bfC =\begin{bmatrix}\bfC_\bfA \\ \bfzero	
\end{bmatrix}\in \mathbb{R}^{m\times n}$ where $\bfC_\bfA \in \bbR^{n \times n}$ is diagonal and
$\bfS\in\mathbb{R}^{p\times n}$ is diagonal (not necessarily
square), satisfying $\bfC^T\bfC+\bfS^T \bfS=\bfI$. The existence of the CS
decomposition can be found in~\cite[Section~22.1]{bjorck1996numerical}. Then the GSVD
of $\{\bfA,\bfL\}$ is defined as follows
\begin{equation} \label{GSVD} \bfA=\bfP \bfC
\underbrace{\bfW^{T}\bfR}_{\bfZ^{-1}} \qquad \mbox{and} \qquad \bfL=\bar{\bfP}\bfS
\underbrace{\bfW^{T}\bfR}_{\bfZ^{-1}},
\end{equation} where for convenience, we have ordered the entries in
$\bfC$ and  $\bfS$ such that
$$1 \geq c_{1}\geq \cdots \geq c_{\min\{n,p\}}\geq 0, \qquad c_{\min\{n,p\}+1} = \cdots = c_n = 1,$$
$$0 \leq s_{1}\leq...\leq s_{\min\{n,p\}}\leq 1.$$ Using the GSVD, the
general-form Tikhonov solution can be written as
 \begin{align}
\bfx_\lambda  &= \sum_{i=1}^n \phi_i
\frac{\bfp_i^T \bfb}{c_i} \bfz_i \label{eq:gsvdfilter}\\
& = \bfZ \bfPhi \begin{bmatrix}
	\bfC_\bfA^{-1} & \bfzero
\end{bmatrix} \bfP\t\bfb \equiv \bfA^\dagger_\lambda \bfb\,, \label{eq:gsvdsoln}
 \end{align}
where $\bfPhi \in \bbR^{n \times n}$ is a diagonal matrix with \emph{filter factors},
\begin{equation}
	\label{eq:Tikfilterfactors}
	\phi_i = \begin{cases} \frac{c_i^2}{c_i^2 + \lambda^2 s_i^2} & \mbox{if } i=1,...,\min\{n,p\} \\
1 &	\mbox{if } i = \min\{n,p\}+1, ..., n
\end{cases}
\end{equation}
on the diagonal.  Notice that for $i \leq \min\{n,p\}$, $\phi_i$ can be written as $\frac{t_i^2}{t_i^2 + \lambda^2}$, where $t_i = \frac{c_i}{s_i}$ are the generalized singular values. It is easy to see that for standard-form Tikhonov regularization where $\bfZ$ and $\bfC$ contain the right singular vectors and singular values of $\bfA$ (with reordering) and $s_i=1$ for all $i$, the filter factors are given by $\phi_i = \frac{c_i^2}{c_i^2 + \lambda^2}$. Therefore, by allowing $\bfL \neq \bfI$, we not only affect the definition of the filter factors, but also alter the basis vectors $\bfz_i$ used to represent the regularized solution $\bfx_\lambda$.

\begin{figure}
	\begin{center}
		\begin{tabular}{cc}
		\includegraphics[scale=0.4]{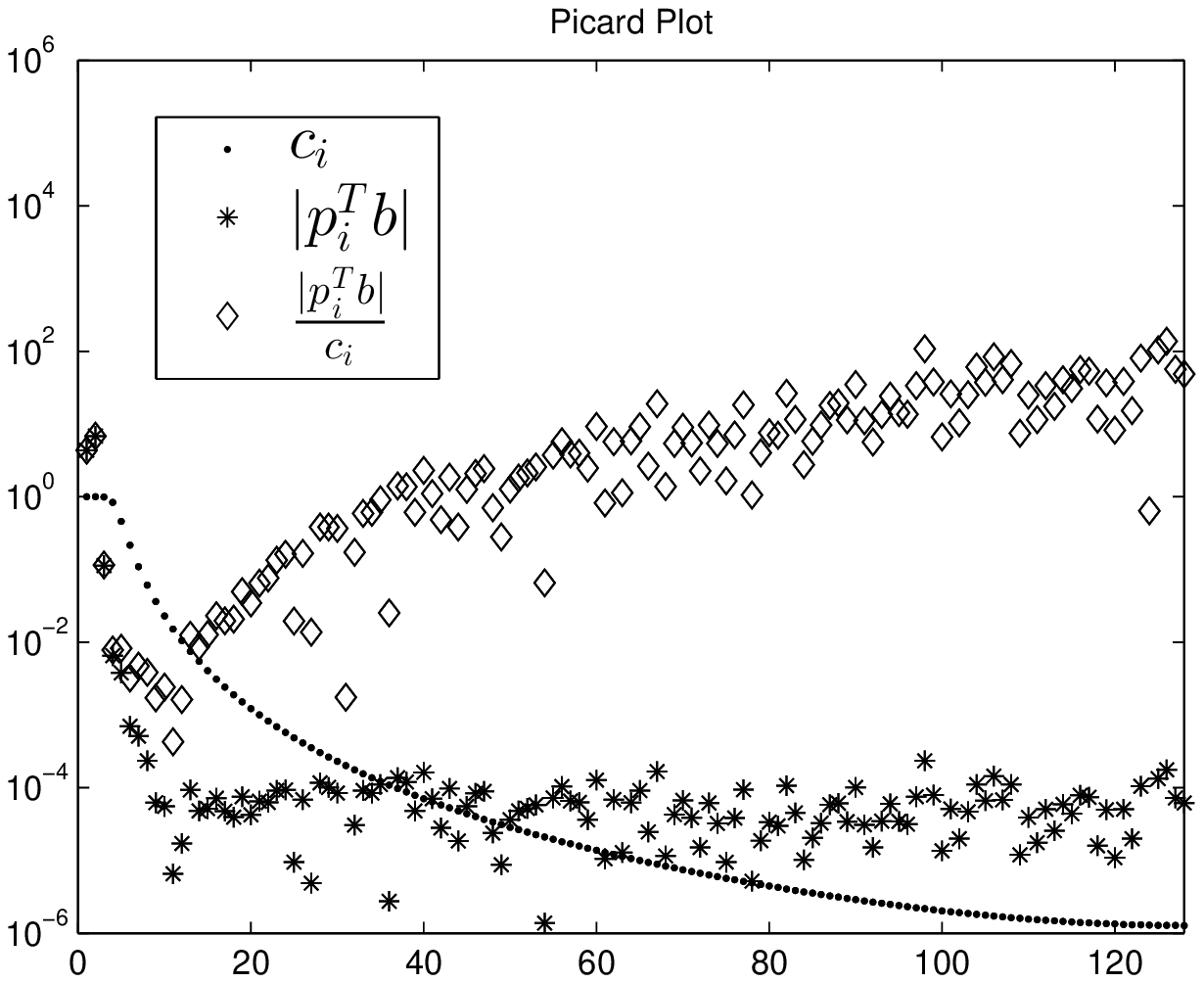} &
        \includegraphics[scale=0.4]{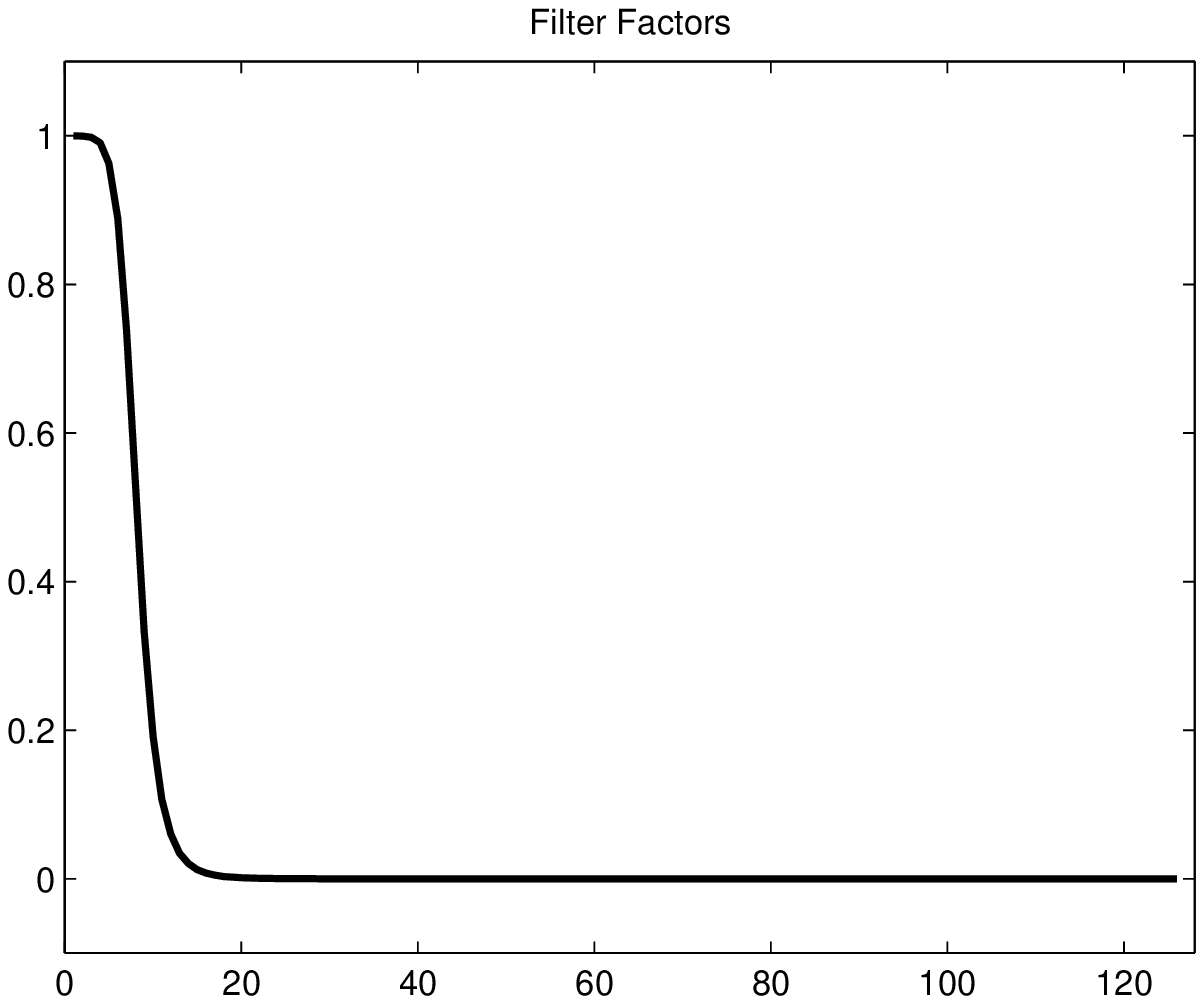}\\
		(a) & (b)
		\end{tabular}
	\end{center}
\caption{A sample Picard plot is provided in (a), where $c_i$ are the generalized singular values and $|\bfp_i ^T \bfb|$ are the absolute values of the spectral coefficients of the observed data. Filter factors are provided in (b) as a function of the index $i$.}
	\label{fig:Picard}	\end{figure}

For ill-posed inverse problems, it is well known that the problem satisfies the discrete Picard condition~\cite{HansenBook1998}.  That is, the values $|\bfp_i^T \bfb|$ decay on average faster than the values $c_i$, until an index is reached where the noise components dominate the solution. After that point, the coefficients $|\bfp_i^T \bfb|$ stabilize around the noise level and $c_i$'s continue decreasing, resulting in amplification of errors in the reconstruction.  A sample Picard plot is provided in Figure~\ref{fig:Picard}. Since the signal
is contained primarily in the subspace spanned by $\bfz_i$ for
small $i$ and including $\bfz_i$ for larger $i$ results in errors, a good value of the regularization parameter $\lambda$ should filter out the terms in~\eqref{eq:gsvdfilter} for larger values of $i$. In other words, the filter factors determined by $\lambda$ should go to zero to counteract the amplification of the error. Figure~\ref{fig:Picard} shows filter factors corresponding to the parameter $\lambda$ that produces the smallest mean squared errors.

Therefore, it is a very important task to find a good regularization parameter $\lambda$. Here we describe two standard methods for selecting regularization parameters for general-form Tikhonov regularization: one that requires an estimate of the level of noise, namely the discrepancy principle, and one that does not require prior knowledge about the noise level, namely GCV.

The basic idea of the discrepancy principle is to select a regularization parameter $\lambda$ such that the residual satisfies
\begin{equation}
	\label{eq:DP}
	\norm[2]{\bfA \bfx_\lambda - \bfb}^2 = \tau \eta, \,
\end{equation}
where $\eta$ is an estimate of the expected value of $\norm[2]{\bfn}^2$, and $\tau>0$ is a user-defined parameter.  For a given $\lambda$, the GSVD (assumed available) can be used to simplify the residual as
\begin{align}
	\norm[2]{\bfA \bfx_\lambda - \bfb}^2 & = 	\norm[2]{\left(\begin{bmatrix}
		\bfPhi &  \nonumber \\ & \bfzero
	\end{bmatrix}-\bfI\right) \bfP\t\bfb}^2 \\
	& =\sum_{i=1}^{\min\{n,p\}} \left(\frac{\lambda^2 s_i^2}{c_i^2 + \lambda^2 s_i^2}\right)^2 (\bfp_i\t \bfb)^2 + \sum_{i={n+1}}^m (\bfp_i\t \bfb)^2\,. \label{eq:residual}
\end{align}
Standard root finding methods can be used to find the regularization parameter $\lambda_{\mbox{DP}}$ that satisfies~\eqref{eq:DP}. However, note that even with a very accurate estimate of $\eta$, results may be sensitive to the choice of $\tau$.

The GCV method \cite{golubGCV} is also a well-known method to compute the regularization parameter $\lambda$. The GCV parameter $\lambda_{\mbox{GCV}}$ is computed to minimize the GCV function,
\begin{equation}
	\label{eq:GCVgeneral}
	\mbox{GCV}(\lambda) = \frac{\norm[2]{(\bfI - \bfA \bfA^\dagger_\lambda)\bfb}^2}{[\trace{(\bfI-\bfA \bfA^\dagger_\lambda)]^2}},
\end{equation}
where $\bfA^\dagger_\lambda$ is defined in~\eqref{eq:gsvdsoln}.  Using the GSVD, the GCV function can be simplified as
\begin{equation}
	\label{eq:GCV}
	\mbox{GCV}(\lambda) = \frac{\displaystyle \sum_{i=1}^{\min\{n,p\}} \left(\frac{\lambda^2 s_i^2}{c_i^2 + \lambda^2 s_i^2}\right)^2 (\bfp_i\t \bfb)^2 + \sum_{i=n+1}^m (\bfp_i\t \bfb)^2}{\displaystyle\left( m-n+ \sum_{i=1}^{\min\{n,p\}} \frac{ \lambda^2 s_i^2}{c_i^2+\lambda^2 s_i^2}\right)^2}\,,
\end{equation}
making evaluation of the GCV function computationally convenient.
As mentioned earlier, computing $\lambda_{\mbox{DP}}$ requires an estimate of the noise level, while computing $\lambda_{\mbox{GCV}}$ does not. In the next section, we describe an approach for computing optimal regularization parameters for cases where training data is available.

\section{Optimal Regularization Parameters}
\label{sec:optfilt}
In this section, we describe a learning approach to compute optimal regularization parameters for the general-form Tikhonov problem~\eqref{GeneralTik} and the multi-parameter Tikhonov problem~\eqref{MultiTik}. We use an empirical Bayes risk minimization framework to incorporate probabilistic information from training data, and we allow optimal regularization parameters to be computed for various error measures.  This approach is an extension of methods developed in~\cite{chung2011designing}.

Suppose we are given training data $(\bfb^{(k)}, \bfx_{\true}^{(k)})$, $k = 1, 2,..., K$, where  $\bfb^{(k)} = \bfA \bfx_{\true}^{(k)}+ \bfn^{(k)}$, and $\bfn^{(k)}$ is unknown. For each observation $\bfb^{(k)}$, we seek either a general-form Tikhonov solution~\eqref{GeneralTik} or a multi-parameter Tikhonov solution~\eqref{MultiTik} and compute the error vector,
\begin{equation}
	\label{eq:errorfunction}
	\bfe^{(k)}(\bflambda) = \bfx_\bflambda^{(k)} - \bfx_{\true}^{(k)}\,,
\end{equation}
where $\bflambda$ is a vector in the multi-parameter case and a scalar in the general-form Tikhonov case.
Let $\rho(\bfxi)\colon\bbR^n \to \bbR_0^+$ be some error measure, for example, $\rho(\bfxi) = ||\bfxi||_p^p$ for $p \geq 1$.  Then the goal is to find the regularization parameter vector $\bflambda$ that minimizes the average error of the reconstructions for the training data,
\begin{equation}
	\label{eqn:optparams}
f_K(\bflambda) \equiv \frac{1}{K} \sum_{k = 1}^{K} \rho(\bfe^{(k)}(\bflambda))\,.
\end{equation}
That is, we seek a solution to the following \emph{empirical} Bayes risk minimization problem,
\begin{equation}\label{eqn:empiricalBayes}
	\min_\bflambda f_K(\bflambda)\,.
\end{equation}

 Notice that underlying~\eqref{eqn:empiricalBayes} is a Bayes risk minimization problem, where the goal is to compute $\bflambda$ that minimizes the expected value of the reconstruction errors over the joint distribution of the noise and unknown parameters \cite{chung2011designing}. Since these distributions are rarely known in practice, we approximate the expected value with the sample mean and instead seek a solution to~\eqref{eqn:empiricalBayes}. Minimizing the average error~\eqref{eqn:optparams} is one approach based on this interpretation that we are minimizing the expected value of the errors, but one could also consider minimizing the median of the errors or the maximum value of the errors, depending on the desired design or goal of the problem.


\subsection{One-parameter General-form Tikhonov} 
\label{sub:one_parameter_tikhonov}
For the one-parameter \linebreak general-form Tikhonov problem, one could use derivative-free methods (e.g., \linebreak \verb|fminsearch| in MATLAB) to solve~\eqref{eqn:empiricalBayes} for various error measures. However, it is also possible to compute the derivative $f_K'(\lambda)$ and solve $f_K'(\lambda)=0$ using standard root finding methods.  For completeness, we provide the derivatives here.

 Let us define vector $\bfgamma \in \bbR^{n \times 1}$ with elements $\gamma_i = \frac{\bfp_i\t \bfb}{c_i}$ and let $\bfGamma\in\bbR^{n \times n}$ be a diagonal matrix with entries $\gamma_i$.  Then decompose
\begin{equation}
	\bfZ = \begin{bmatrix} \widehat \bfZ, \widetilde \bfZ
	\end{bmatrix}, \quad \bfgamma = \begin{bmatrix} \widehat \bfgamma \\ \widetilde \bfgamma	\end{bmatrix}, \quad \mbox{and} \quad \bfGamma = \begin{bmatrix} \widehat \bfGamma & \\ & \widetilde \bfGamma
\end{bmatrix}\,,
\end{equation} where $\widehat \bfZ \in \bbR^{n \times \min\{n,p\}}$, $\widehat \bfgamma \in \bbR^{\min\{n,p\} \times 1}$, $\widehat \bfGamma\in \bbR^{\min\{n,p\} \times  \min\{n,p\}}$ and $\widetilde \bfZ, \widetilde \bfgamma,$ and $\widetilde \bfGamma$ have corresponding size.  Then the solution to the general-form Tikhonov problem can be written as
\begin{equation}
	\label{eq:gsvdsoln2}
\bfx_\lambda^{(k)}  = \sum_{i=1}^{n} \phi_i
\frac{\bfp_i^T \bfb^{(k)}}{c_i} \bfz_i  = \widehat \bfZ \widehat \bfGamma \bfphi(\lambda)+ \widetilde \bfZ \widetilde \bfgamma,
\end{equation}
where $\bfphi \in \bbR^{\min\{n,p\} \times 1}$ contains filter factors $\bfphi_i$ which depend on $\lambda$. For each $k,$ the Jacobian of the errors with respect to $\lambda$ can be written as
\begin{equation} \label{eq:jac}
\bfJ^{(k)} = \widehat \bfZ \widehat \bfGamma^{(k)} \bfpsi\,,
\end{equation}
where the $i$th element of vector $\bfpsi$ is the derivative of $\phi_{i}$ with respect to $\lambda,$ which is given by
\begin{equation}
	\psi_i =  \frac{-2 \lambda c_i^2 s_i^2}{(c_i^2 +  \lambda^2 s_{i}^2)^2}.
\end{equation}
Thus, the derivative $f_K'(\lambda)$ can be computed as
\begin{equation}
	f_K'(\lambda) = \frac{1}{K} \sum_{k=1}^K {\bfJ^{(k)}}\t \nabla_{\rho}(\bfe^{(k)}) \,,
\end{equation}
where the $\nabla_{\rho}$ contains the partial derivatives of $\rho$ evaluated at $\bfe^{(k)}(\lambda).$  We remark that the second derivative, if desired, could be numerically approximated.

It is worthwhile to note that one could compute \emph{optimal error} filters as described in \cite{chung2011designing}, where no functional form for the filter factors is assumed.  In this case, the basis vectors for the solution would be the columns of $\bfZ$, and $n$ filter factors would need to be computed rather than one regularization parameter, thereby requiring not only more computational cost for the minimization but also more training data.

\subsection{Multi-parameter Tikhonov} 
\label{sub:multiple_parameter_tikhonov}

In this section we describe an approach to compute optimal or near-optimal regularization parameters for the multi-parameter Tikhonov problem~\eqref{MultiTik}, whose solution can be written as,
\begin{equation}
	\label{eq:multiTikNormal}
  \bfx_{\bflambda} =
  \left(\bfA^T\bfA + \sum_{j=1}^J \lambda_j \bfL^T_j\bfL_j\right)^{-1}\bfA^T\bfb.
 \end{equation}

Unlike the one-parameter general-form Tikhonov problem where the GSVD provides a decomposition of $\{\bfA, \bfL\}$ that allows the solution to be written as a filtered solution~\eqref{eq:gsvdfilter}, a factorization for the multi-parameter Tikhonov problem where all matrices share the same right factor,
\begin{equation}
		\label{eq:simultdiag}
	\bfA  = \bfP \bfC \bfZ^{-1} \quad \mbox{and} \quad	\bfL_j  = \bar\bfP_j \bfS_j \bfZ^{-1} \quad j=1,...,J \,,
	\end{equation}
may not be possible, in general \cite{BrezinskiMulti2003}.  Thus, a filtered representation for $\bfx_\bflambda$ is not always possible.

For problems where such a decomposition does not exist, a bilevel optimization approach could be used to compute optimal parameters~\cite{colson2007overview,kunisch2013bilevel}.
However, this approach would require computing \eqref{eq:multiTikNormal} for various estimates of $\bflambda$, which could become computationally expensive. Another approach that was proposed in~\cite{ChungKilmerOleary} uses operator approximations to estimate regularization parameters. For problems where operator approximations can be found that satisfy the above simultaneous diagonalizability condition~\eqref{eq:simultdiag}, we propose to use the approximate problem to estimate near-optimal regularization parameters and then solve the original problem using the computed parameters. Thus, here we restrict our derivation to square matrices, i.e., $m=n$ and $p_j=n$ for $j=1,...,J,$ and consider scenarios where such a decomposition exists or can be approximated.


We remark that in signal and image processing applications, it is common to encounter matrices $\bfA$ and $\bfL_j$ that satisfy~\eqref{eq:simultdiag}, under certain assumptions on the boundary conditions.
In these situations, $\bfP, \bar\bfP_j$ and $\bfZ^{-1}$ typically represent the same frequency transform, and we get the following simultaneous diagonalizability condition
\begin{equation}
	\label{eq:simultdiag2}
\bfA = \bfQ \bfC \bfQ^* \quad \mbox{and} \quad \bfL_j = \bfQ \bfS_j \bfQ^*\quad j=1,...,J \,,
\end{equation}
where $\bfQ$ is an orthogonal (or unitary) matrix.
For instance, in image deblurring, if we assume \emph{periodic boundary conditions} on the image, then the blurring operator $\bfA$ and the regularization matrices $\bfL_j$ are block circulant matrices with circulant blocks (BCCB). In this case, $\bfQ^*$ represents the 2D discrete Fourier transform (DFT) matrix. If we assume \emph{reflexive boundary conditions} and the point spread functions defining $\bfA$ and $\bfL_j$ satisfy a double symmetry condition \cite{DeblurringBook}, then the matrices can each be written as
a sum of BTTB (block Toeplitz with Toeplitz blocks), BTHB (block Toeplitz with Hankel blocks), BHTB (block Hankel with Toeplitz blocks), and BHHB (block Hankel with Hankel blocks) matrices. In this case, the matrices are diagonalized by the  discrete cosine transform (DCT), where $\bfQ^*$ represents the 2D discrete DCT matrix.

Next, we show that with the above decomposition~\eqref{eq:simultdiag2}, the multi-parameter Tikhonov solution can be written as a filtered solution.  For notational purposes, $|\bfDelta|^2$ represents element-wise absolute value followed by element-wise square for any matrix $\bfDelta$.  Then from~\eqref{eq:multiTikNormal}, the solution of~\eqref{MultiTik} can be written as,
\begin{align}
  \bfx_{\bflambda} &=
  \bfQ \left(|\bfC|^2\left (|\bfC|^2  + \sum_{j=1}^J \lambda_j^2|\bfS_j|^2\right)^{-1}\right)\bfC^{-1}\bfQ^* \bfb \nonumber \\
& = \sum_{i=1}^{n} \phi_i
\frac{\bfq_i^* \bfb}{c_i} \bfq_i \nonumber\\
& = \bfQ \bfPhi \bfC^{-1}\bfQ^* \bfb \,, \label{eq:filteredsolnmulti}
 \end{align}
where as before $\bfPhi$ is a diagonal matrix containing the filter factors, which in this case are given by
\begin{equation}
	\label{eq:filterfactorsmulti}
	\phi_i = \frac{|c_i|^2}{|c_i|^2 +  \sum_{j=1}^J \lambda_j^2 |s_{i,j}|^2},
\end{equation}
with $s_{i,j}$ being the $i$th diagonal element of matrix $\bfS_j$.

In this case, the GCV method can provide a vector of regularization parameters $\bflambda_{\mbox{GCV}}=[\lambda_1,\dots, \lambda_J]$ by minimizing the GCV function \eqref{eq:GCVgeneral} where now $\bfA^\dagger_\lambda=\bfQ \bfPhi \bfC^{-1}\bfQ^*$. Using decomposition~\eqref{eq:simultdiag2}, the GCV function can be written as
\begin{equation}
	\label{eq:GCVmultiGSVD}
	\mbox{GCV}(\bflambda) = \frac{\displaystyle \sum_{i=1}^{n} \left(1-\phi_i\right)^2 (\bfq_i^* \bfb)^2}{\displaystyle\left(\sum_{i=1}^{n} 1-\phi_i\right)^2}\,.
\end{equation}
In the numerical results section, we will use the GCV regularization parameters for comparison.

For notational convenience, let $\bfGamma$ be an $n \times n$ diagonal matrix whose $i$th diagonal element is $\gamma_i = \frac{\bfq_i^* \bfb}{c_i}$.  Then~\eqref{eq:filteredsolnmulti} can be written as
 \begin{equation}
\bfx_\bflambda = \bfQ \bfGamma \bfphi(\bflambda),
 \end{equation}
where $\bfphi(\bflambda) \in \bbR^{n\times 1}$ is the vector of filter factors~\eqref{eq:filterfactorsmulti} that depends on the vector of regularization parameters $\bflambda$.

Recall that our goal is to compute optimal regularization parameters by solving~\eqref{eqn:empiricalBayes}, where instead of one regularization parameter $\lambda$, we now have multiple parameters $\lambda_1,\dots, \lambda_J$, and therefore the error in~\eqref{eqn:optparams} is given by
\begin{equation}
	\bfe^{(k)}(\bflambda) = \bfQ \bfGamma^{(k)} \bfphi(\bflambda) - \bfx_{\true}^{(k)}\,.
\end{equation}
Here, for each training data $\bfb^{(k),}$ we define $\bfGamma^{(k)}$ where the $i$th diagonal element of $\bfGamma^{(k)}$ is $\gamma_i^{(k)}=\frac{\bfq_i^* \bfb^{(k)}}{c_i}$.

Since we must optimize over several parameters, we propose to use a Gauss-Newton approach to solve~\eqref{eqn:empiricalBayes}.  The Jacobian of the errors with respect to $\bflambda$ can be written as
\begin{equation} \label{eq:multjac}
\bfJ = \frac{1}{K}\left[\begin{array}{c}\bfJ^{(1)} \\ \vdots \\ \bfJ^{(K)} \end{array}\right] \,, \quad \mbox{ with }  \bfJ^{(k)} = \bfQ \bfGamma^{(k)} \bfPsi\,,
\end{equation}
where the $(i,j)$-th entry of the matrix $\bfPsi \in \bbR^{n\times J}$ is the derivative of $\phi_{i}$ with respect to $\lambda_j.$
Define $\bar\bfS \in \bbR^{n \times J}$ whose $i,j$ entry is $|s_{i,j}|^2$, then the matrix $\bfPsi$ is given by,
\begin{equation}
	\bfPsi = -2  \bfT \, \bar\bfS \, \bfLambda\,,
\end{equation}
where $\bfLambda = \diag{\bflambda}$
and $\bfT$ is a diagonal matrix whose $i$th diagonal entry is given by $\frac{|c_i|^2}{(|c_i|^2 +  \sum_{j=1}^J \lambda_j^2 |s_{i,j}|^2)^2}$.
The gradient of $f_K(\bflambda)$
and the Gauss-Newton approximation of the Hessian can be written as
\begin{equation}
	\label{eq:multGNHessian}
	\bfg = \frac{1}{K} \sum_{k=1}^K {\bfJ^{(k)}}^* \nabla_{\rho}(\bfe^{(k)}) \,, \quad \mbox{and} \quad \bfH = \frac{1}{K} \sum_{k=1}^K {\bfJ^{(k)}}^* \bfF^{(k)} \bfJ^{(k)} \,,
\end{equation}
where the $\nabla_{\rho}(\bfe^{(k)})$ contains the partial derivatives of $\rho$ evaluated at $\bfe^{(k)}(\lambda),$ and $\bfF^{(k)}$ contains the second derivatives (i.e., the Hessian) of $\rho$ evaluated at $e^{(k)}(\bflambda).$

For the special case where $\rho(\bfxi) = \frac{1}{2}\norm[2]{\bfxi}^2,$ we have $\nabla_{\rho}(\bfe^{(k)}) = \bfe^{(k)}$ and $\bfF^{(k)}=\bfI$, so the Hessian approximation simplifies nicely,
\begin{align*}
\bfH & = \frac{1}{K} \sum_{k=1}^K {\bfJ^{(k)}}^*  \bfJ^{(k)}  = \frac{1}{K} \sum_{k=1}^K \bfPsi^* \left(\bfGamma^{(k)}\right)^* \bfQ^* \bfQ   \bfGamma^{(k)} \bfPsi  \\
& = \frac{4}{K} \bfLambda \bar \bfS^* \left\{ \bfT \left(\sum_{k=1}^K |\bfGamma^{(k)}|^2\right) \bfT\right\} \bar\bfS \bfLambda \,,
\end{align*}
where we note that the matrix contained in the curly brackets is diagonal.
Notice that for the case $J=1$, these derivations are equivalent to those for the one-parameter general-form Tikhonov case.


\section{Numerical Results}
\label{sec:numerics}
In this section, we provide some numerical results that illustrate the performance of the optimal regularization parameters for general-form Tikhonov and multi-parameter Tikhonov regularization.  We consider three investigations.  Example 1 is a 1D signal deconvolution example that compares optimal Tikhonov and optimal error filters for both standard SVD and generalized SVD reconstructions.  Examples 2 and 3 are 2D image deconvolution examples.  Example 2 compares one $\bfL$ versus multiple $\bfL$s for various error measures, and Example 3 investigates an example where approximate matrices $\bfA$ and $\bfL$ can be used to obtain regularization parameters for the original problem.

\paragraph{Example 1} In this first example, we use a 1D signal deconvolution example to investigate the benefits of incorporating a regularization matrix $\bfL$ in the optimal filter framework, compared to the standard Tikhonov case where $\bfL=\bfI$. We compare optimal Tikhonov filters and optimal error filters for each case.

\begin{figure}
	 \caption{Columns of these images were used as training signals for computing optimal regularization parameters for Example 1.}
 	\label{fig:OneLtraining}
 	\begin{center}
 		\includegraphics[scale=0.3]{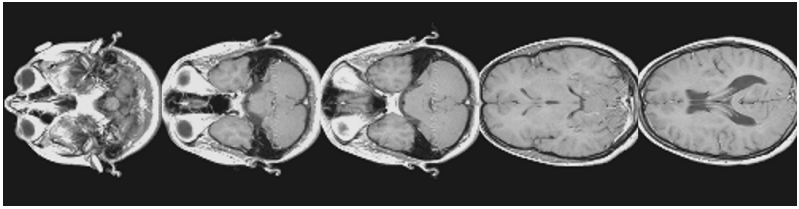}
 	\end{center}		
 \end{figure}

To generate training signals, we took $200$ columns of each of the $5$ MRI images shown in Figure~\ref{fig:OneLtraining}, thereby resulting in $1,000$ training signals, each of size $256 \times 1$.  Each signal was blurred with a Gaussian point spread function with mean 0 and variance 1.  Gaussian white noise was added, where noise levels were randomly selected between $0.2$ and $0.25$.  That is, a noise level of $0.2$ means $\norm[2]{\bfeta}^2/\norm[2]{\bfA \bfx_\true}^2 = 0.2$. Validation signals were generated in the same manner using columns from different MRI images.

For this example, we select regularization matrix $\bfL\in \bbR^{257 \times 256}$ to represent a discretization of a first derivative operator,
\begin{equation}L=
\left(\begin{array}{cccc}
     1 &  & &  \\
     -1 & 1 & &  \\
       & \ddots & \ddots  &   \\
      &   & -1 & 1 \\
      & &   & -1 \\
  \end{array}\right).
\end{equation}

We compare four approaches that use training data to obtain optimal filters.
\begin{itemize}
	\item opt-Tik-SVD and opt-Tik-GSVD correspond to standard and general \\Tikhonov solutions obtained using filter factors \eqref{eq:Tikfilterfactors}, where $\lambda$ minimizes the average mean squared error over all training data, i.e., $\lambda$ solves~\eqref{eqn:empiricalBayes} with $\rho(\bfxi) = \norm[2]{\bfxi}^2$.
	\item opt-error-SVD and opt-error-GSVD correspond to solutions that use the optimal error filters for the SVD and GSVD bases respectively.  That is, the optimal error filter factors assume no functional form, so optimization is performed over $256$ independent filter factors.
\end{itemize}
It is worth noting that opt-Tik-SVD and opt-error-SVD filters were introduced and compared in \cite{chung2011designing}.  Our goal here is to investigate their GSVD counterparts.

\begin{figure}[bthp]
	\caption{Box and whisker plots for the 2-norm of the reconstruction errors for the validation set in Example 1.  Horizontal lines in the box correspond to the median and 25th and 75th percentiles.  Whiskers correspond to extreme data points, and outliers are plotted individually.}
	\label{fig:OneLbox}	
	\begin{center}
		\includegraphics[scale=0.3]{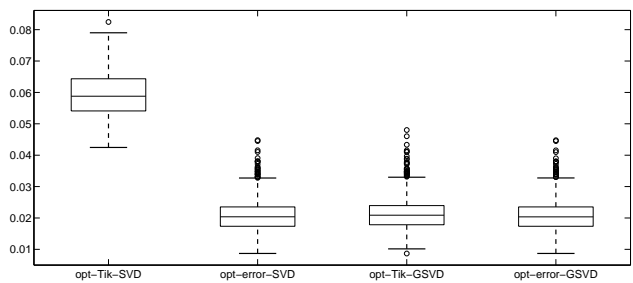}
	\end{center}	
\end{figure}

Once computed, the optimal filters were used to reconstruct each of the validation signals.  The distribution of the relative errors, $\norm[2]{\bfx_\lambda - \bfx_\true}^2/\norm[2]{\bfx_\true}^2$, for the validation set is presented in Figure~\ref{fig:OneLbox} using box and whisker plots.  The box part presents the median value, along with the 25th and 75th percentiles.  The whiskers correspond to extreme data points and the outliers are plotted individually.  As observed in \cite{chung2011designing}, relative errors for opt-Tik-SVD are consistently higher than those for opt-error-SVD due to the high noise level in the problem.  The optimal Tikhonov filter corresponding to GSVD can produce errors that are comparable to both opt-error-SVD and opt-error-GSVD.  This is an important point because both opt-error-SVD and opt-error-GSVD require $n$ unknowns, whereas opt-Tik-GSVD only requires one.  Of course, the performance of opt-Tik-GSVD relies on a good choice of $\bfL$; such concerns will be addressed in later examples of the multi-parameter case.

\begin{figure}[bthp]
	\caption{Computed optimal filter factors for approaches compared in Example 1.}
 	\label{fig:OneLfilterfactors}	
 	\begin{center}
 		\includegraphics[scale=0.4]{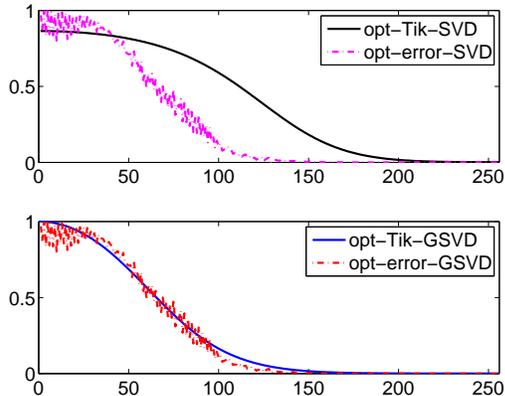}
 	\end{center} 	
 \end{figure}

Filter factors for the four approaches are presented in Figure~\ref{fig:OneLfilterfactors}.
In the top plot, we observe that opt-Tik-SVD has difficulty in obtaining filter factors close to 1 (corresponding to large singular values).  Since Tikhonov filter factors depend on one regularization parameter, there is an inherent limitation of standard-form Tikhonov for problems with high noise levels.  However, this is not true of Tikhonov filter factors with GSVD, as evident in the bottom plot in Figure~\ref{fig:OneLfilterfactors}.  Thus, one of the potential benefits of using general-form Tikhonov is being able to obtain filter factors that can overcome the limitations of standard-form Tikhonov for problems with high noise level.

\begin{figure}[bthp]
	\caption{Absolute error images between the reconstructed signals and the true signals for a subset of the validation signals in Example 1. Images are displayed in inverted colormap so that white corresponds to zero error.}
 	\label{fig:OneLerrorimages}	
 	\begin{center}
 		\includegraphics[scale=0.6]{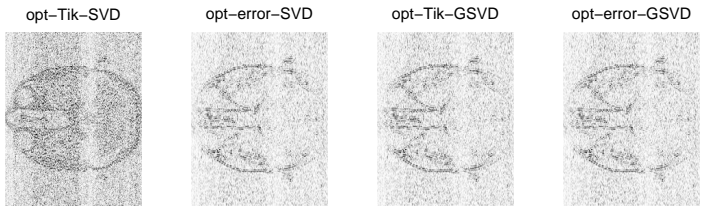}
 	\end{center} 	
 \end{figure}

For a subset of the validation signals, we provide absolute error images between the computed and true solutions in Figure~\ref{fig:OneLerrorimages} in inverted colormap.  White regions in the image correspond to low absolute errors and darker regions correspond to larger errors.  These results are consistent with our observations from Figure~\ref{fig:OneLbox}.

We remark that the training set can be used to obtain uncertainty estimates for the reconstructions.  The average reconstruction errors with corresponding standard deviations for both the training and validation sets are provided in Table~\ref{tab:OneLerrors}.  For comparison purposes, we also include results for 
MSE-Tik-GSVD, which corresponds to general-form Tikhonov, where the filter has the form~\eqref{eq:Tikfilterfactors}, and regularization parameter $\lambda$ minimizes the mean squared error between the true signal and the reconstruction, $\norm[2]{\bfx_\lambda - \bfx_\true}^2$.  This approach requires the true solution from the validation set and is not practical in applications.
Since this parameter depends on the individual problem, these results had to be computed separately for each training and validation signal. From Table~\ref{tab:OneLerrors}, we see that errors corresponding to the training set provide fairly accurate estimates for the validation set.  We observe that opt-error-GSVD produces slightly better results than opt-error-SVD, implying that it is not just the filter factors, but also the basis vectors that play a role in improving the reconstructions.


 \begin{table}
	\caption{Average relative reconstruction error (RRE) and standard deviation (STD) for the training and validation sets in Example 1.}
 	\label{tab:OneLerrors}	
 	\begin{center}
 \begin{tabular}{|l||c|c||c|c|} \hline
 	&  \multicolumn{2}{c||}{Training Set}   & \multicolumn{2}{|c|}{Validation Set}\\ \hline
 	              &    average RRE      &  STD   &     average RRE      &  STD    \\ \hline
 	opt-Tik-SVD  &5.912e-02  & 6.713e-03 & 5.929e-02  & 6.839e-03 \\ \hline
 	opt-error-SVD & 2.145e-02 &5.931e-03  & 2.092e-02 & 5.182e-03 \\ \hline
 	opt-Tik-GSVD    & 2.202e-02  & 5.986e-03 & 2.143e-02 & 5.277e-03\\ \hline
 	opt-error-GSVD   &2.144e-02  &5.927e-03 & 2.091e-02 & 5.180e-03\\ \hline
 	MSE-Tik-GSVD & 2.142e-02  &5.185e-03 & 2.091e-02 &  4.855e-03 \\ \hline
 	\end{tabular}
 	\end{center}
  \end{table}

\begin{figure}[bthp]
	\caption{Pareto curves that display the average relative reconstruction errors for the validation set as a function of the number of training signals used to compute the optimal filters in Example 1.}
 	\label{fig:pareto}	
 	\begin{center}
 		\includegraphics[scale=0.4]{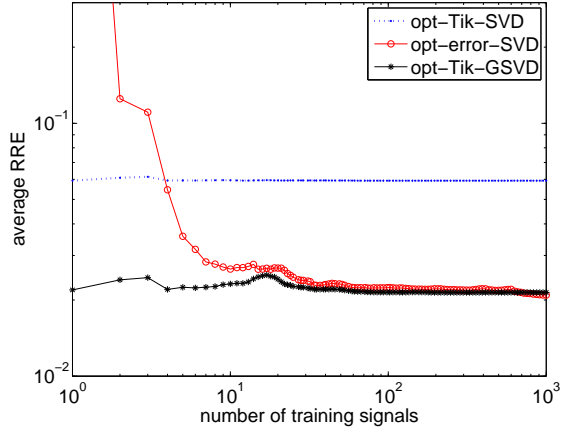}
 	\end{center} 	
 \end{figure}

Lastly, for this example, we investigate the reconstruction errors for the validation set as a function of the number of training signals used to obtain the optimal filters.  We compare opt-Tik-SVD, opt-error-SVD, and opt-Tik-GSVD, and provide average reconstruction errors for the validation set in Figure~\ref{fig:pareto}.  These plots are often referred to as Pareto curves.  We see that average relative reconstruction errors for opt-Tik-SVD and opt-Tik-GSVD remain fairly stable as the number of training signals increases.  On the other hand, for small numbers of training signals, opt-error-SVD is significantly biased towards the training data, but reconstruction errors plateau with enough training signals.  It is interesting to note that opt-error-SVD requires $298$ training data to achieve the same relative reconstruction error that opt-Tik-GSVD can obtain with only $1$ training signal. Thus, opt-Tik-GSVD can be a good alternative to opt-error-SVD if a good $\bfL$ is known a priori and if there are very few signals in the training set.

 \paragraph{Example 2} In this example, we consider a 2D image deblurring example and investigate optimal regularization parameters for Tikhonov regularization, for both the one parameter and multiple parameter case, and we compare different choices of error measures $\rho$.

\begin{figure}
	\caption{Sample training and validation images used in Example 2 can be found in the top and bottom rows respectively.}
 	\label{fig:images}	
 	\begin{center}
		\begin{tabular}{cccc}
	\includegraphics[scale=0.2]{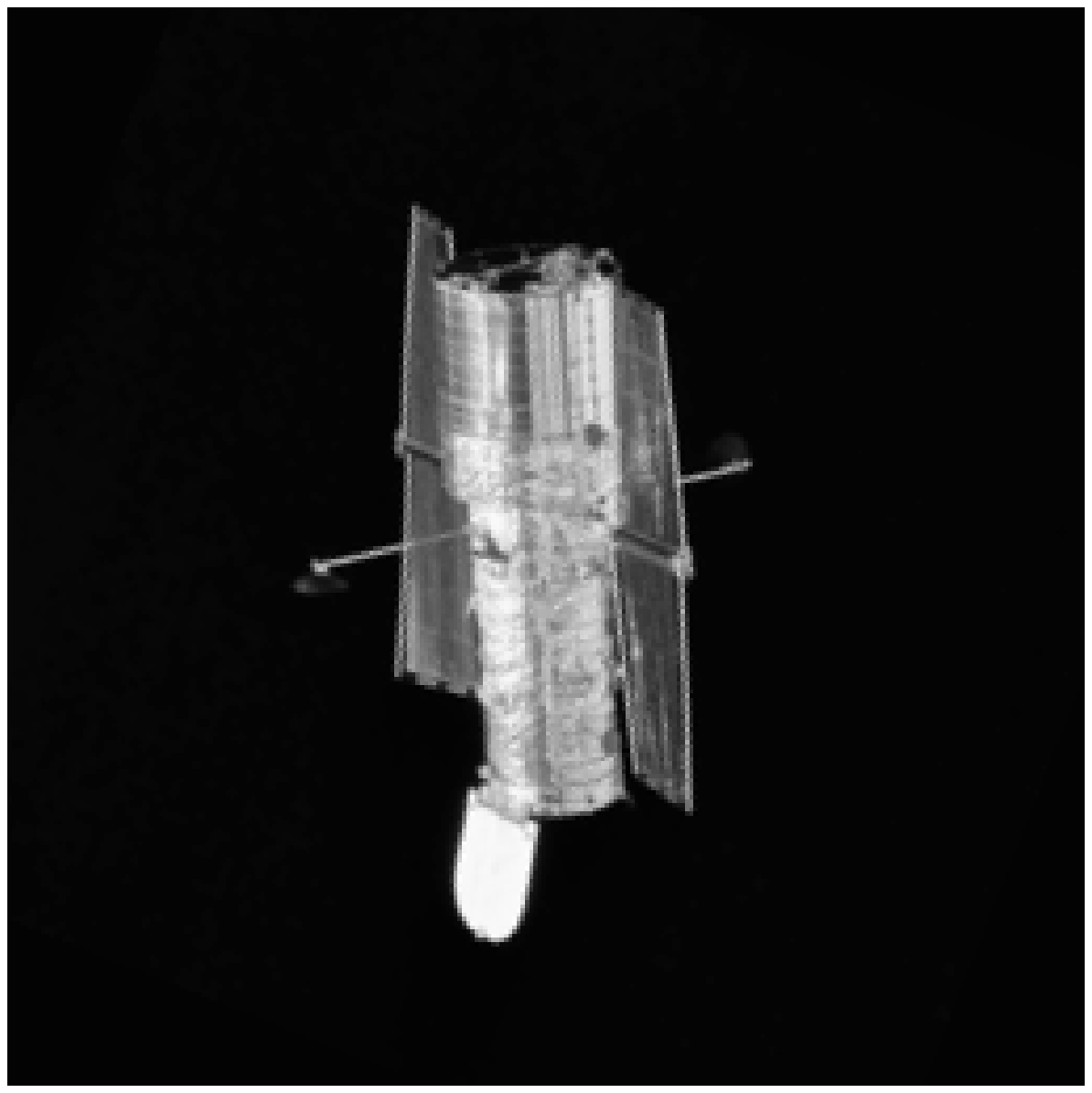} &
 		\includegraphics[scale=0.2]{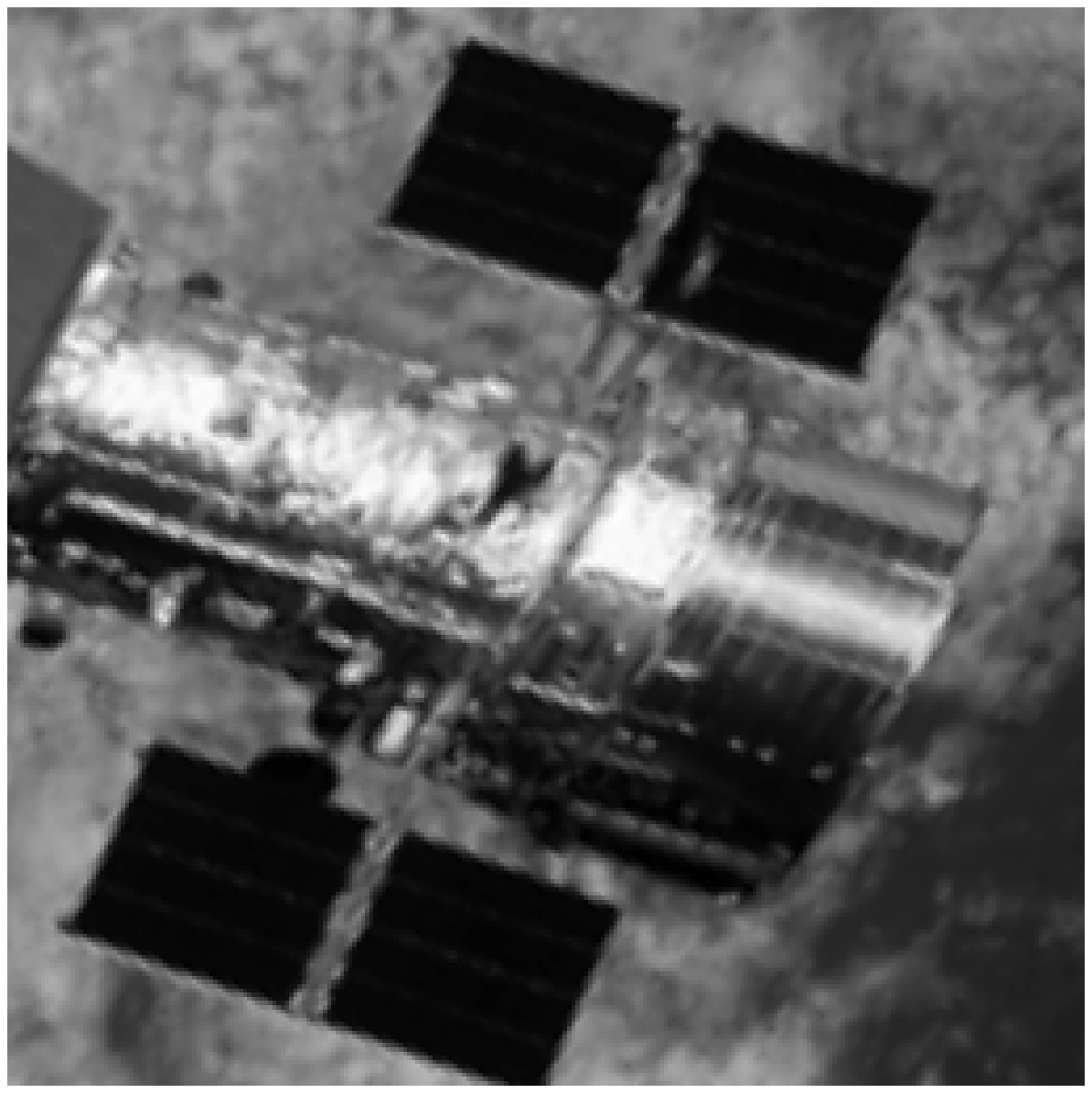} &
 		\includegraphics[scale=0.2]{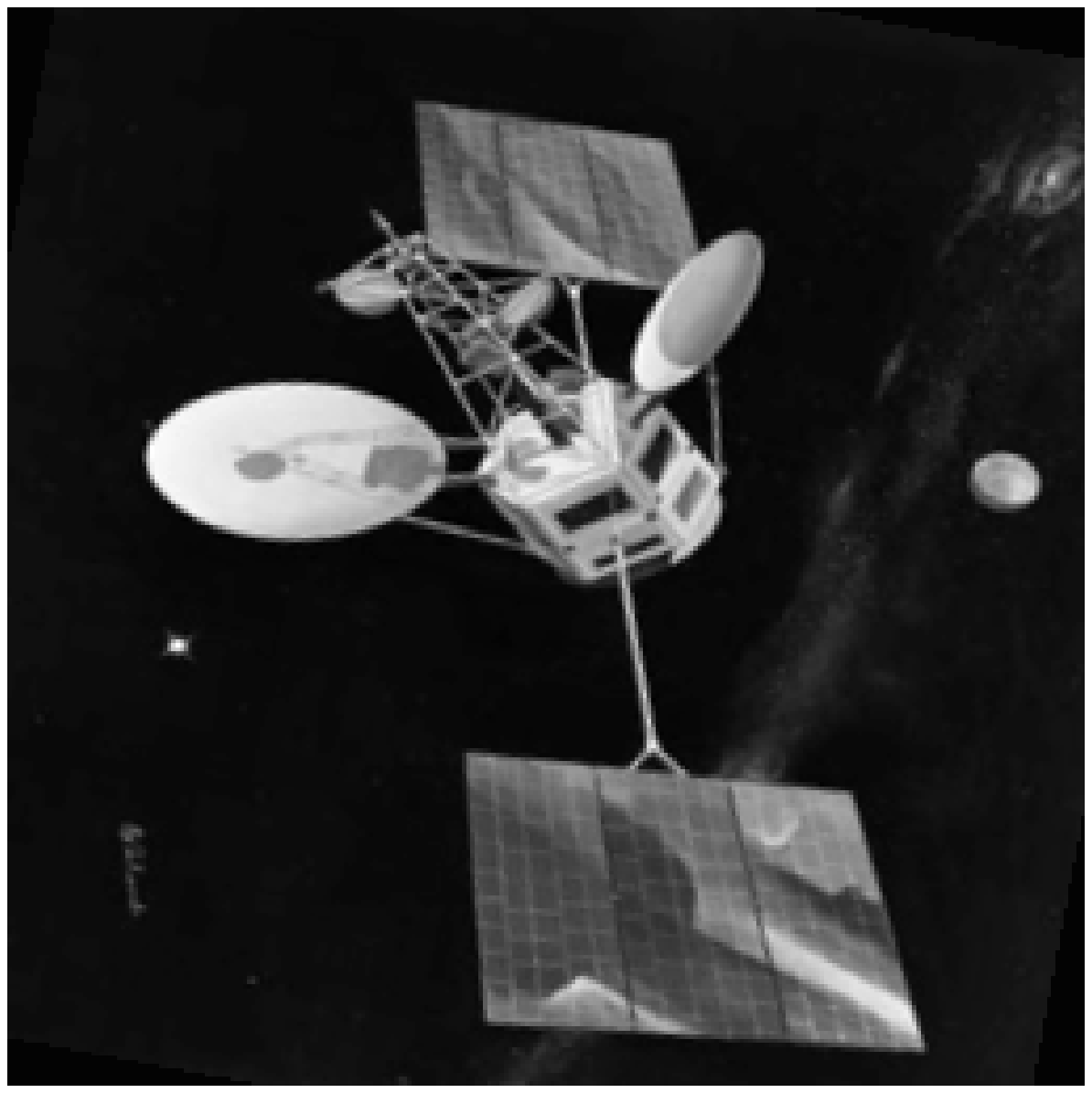} &
 		\includegraphics[scale=0.2]{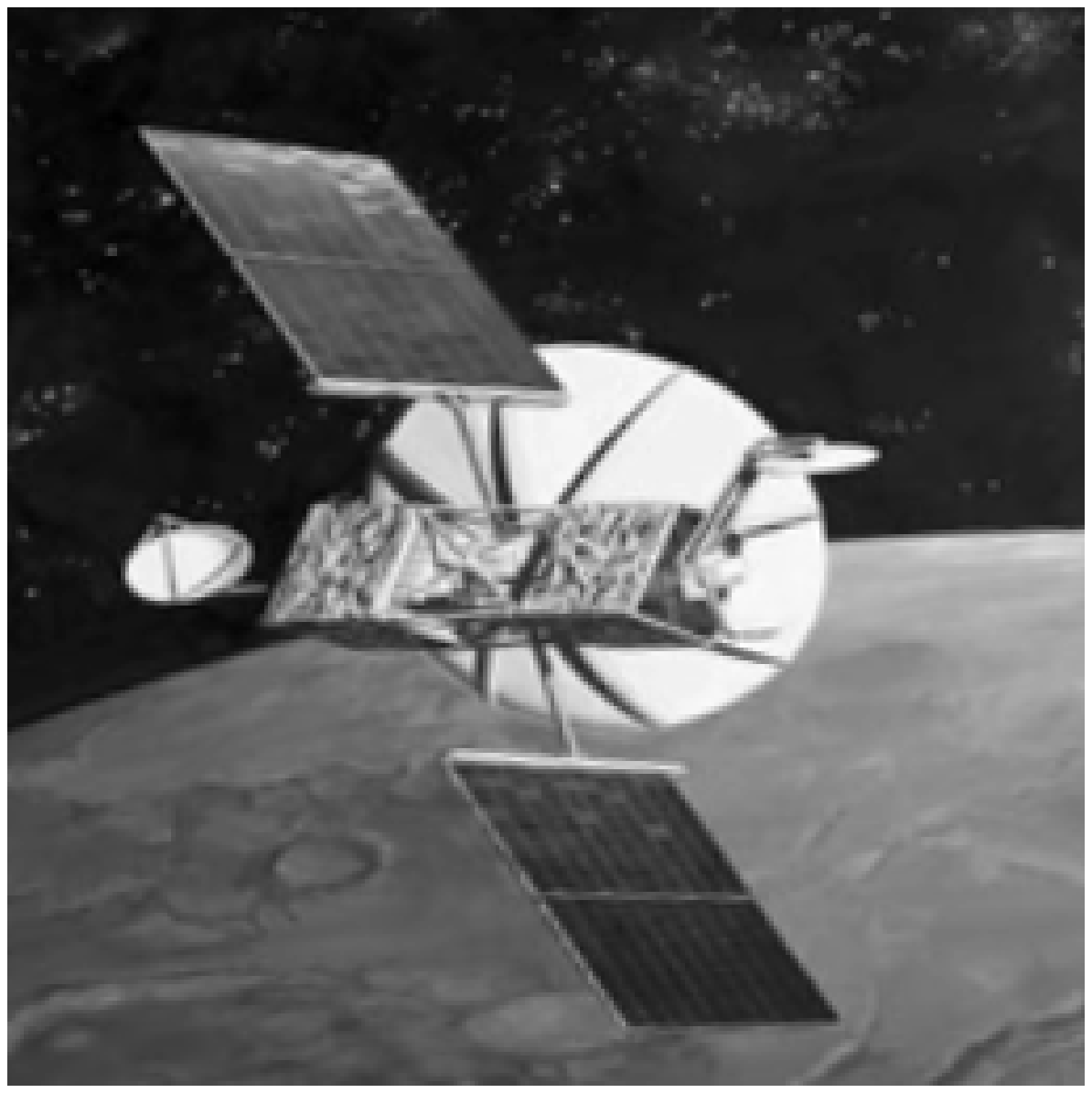} \\
\includegraphics[scale=0.2]{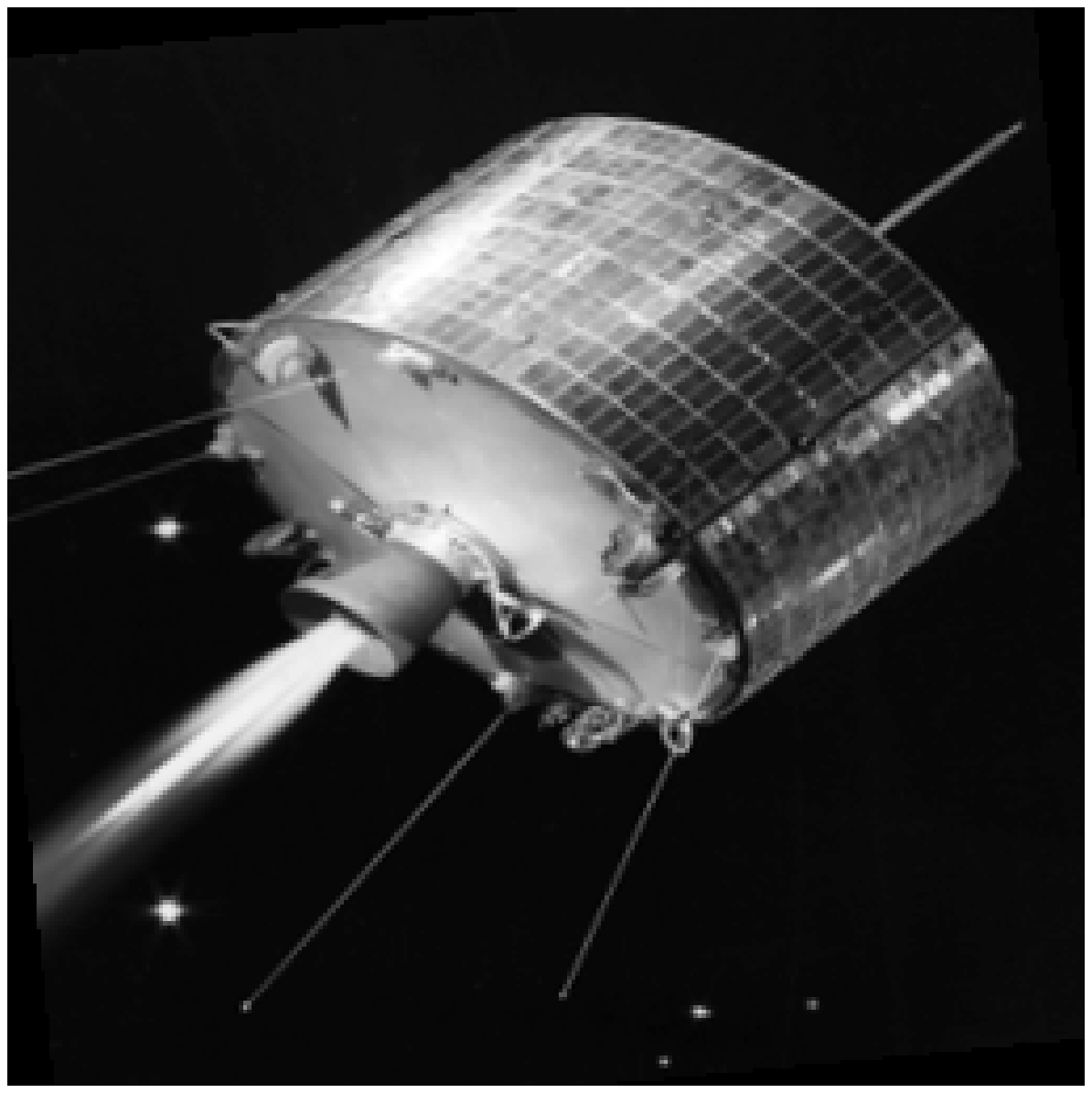} &
 		\includegraphics[scale=0.2]{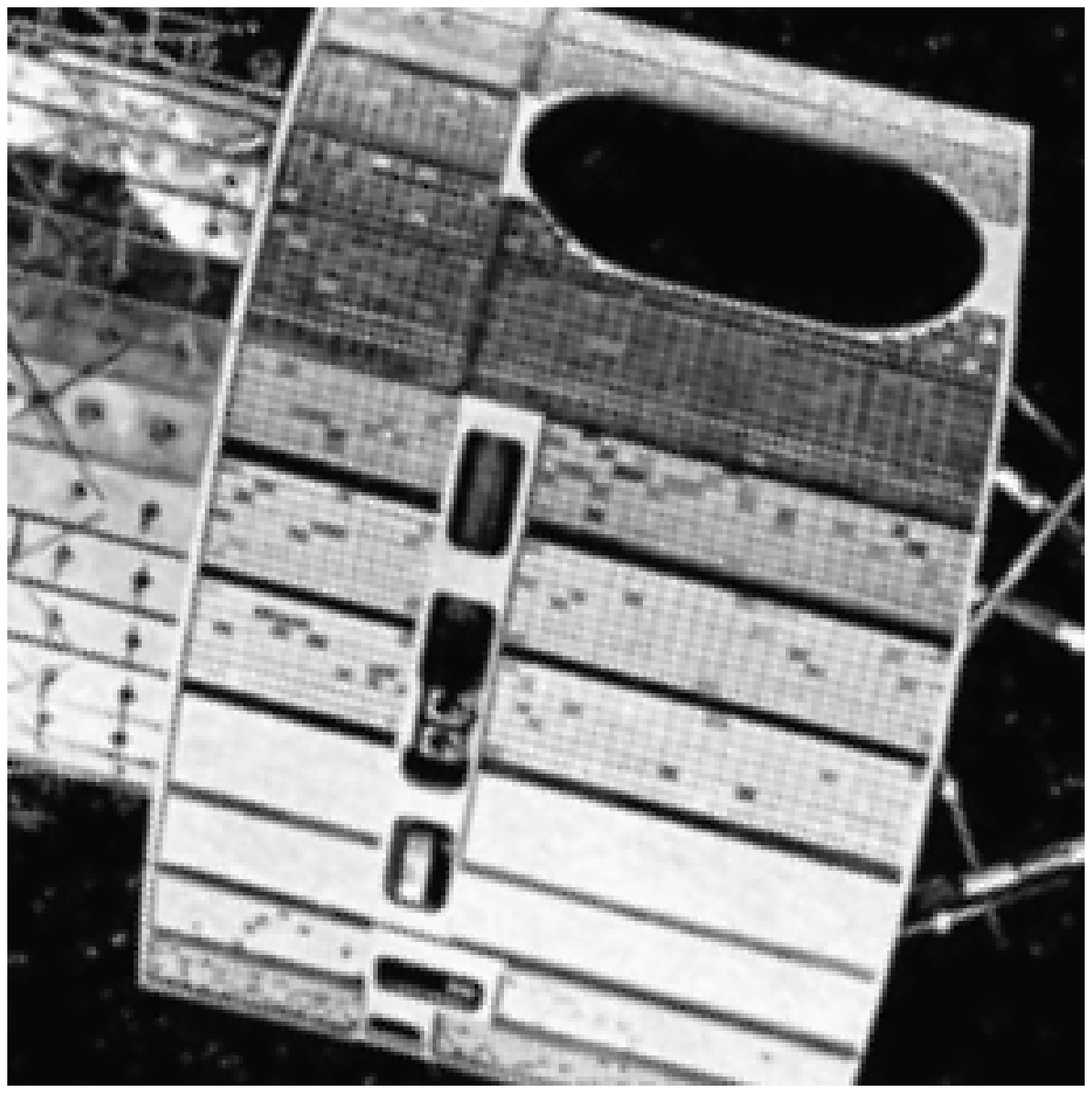} &
 		\includegraphics[scale=0.2]{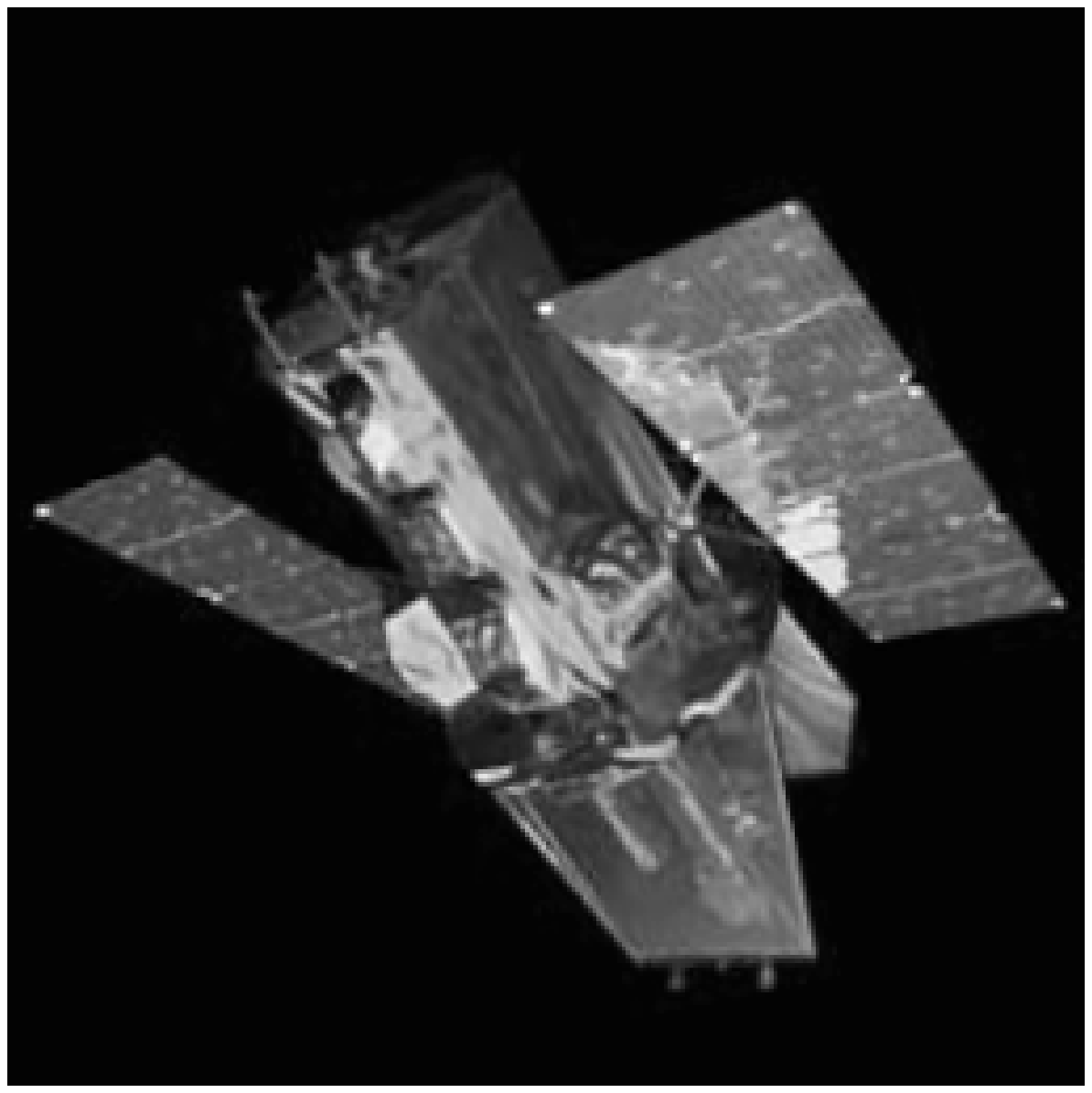} &
 		\includegraphics[scale=0.2]{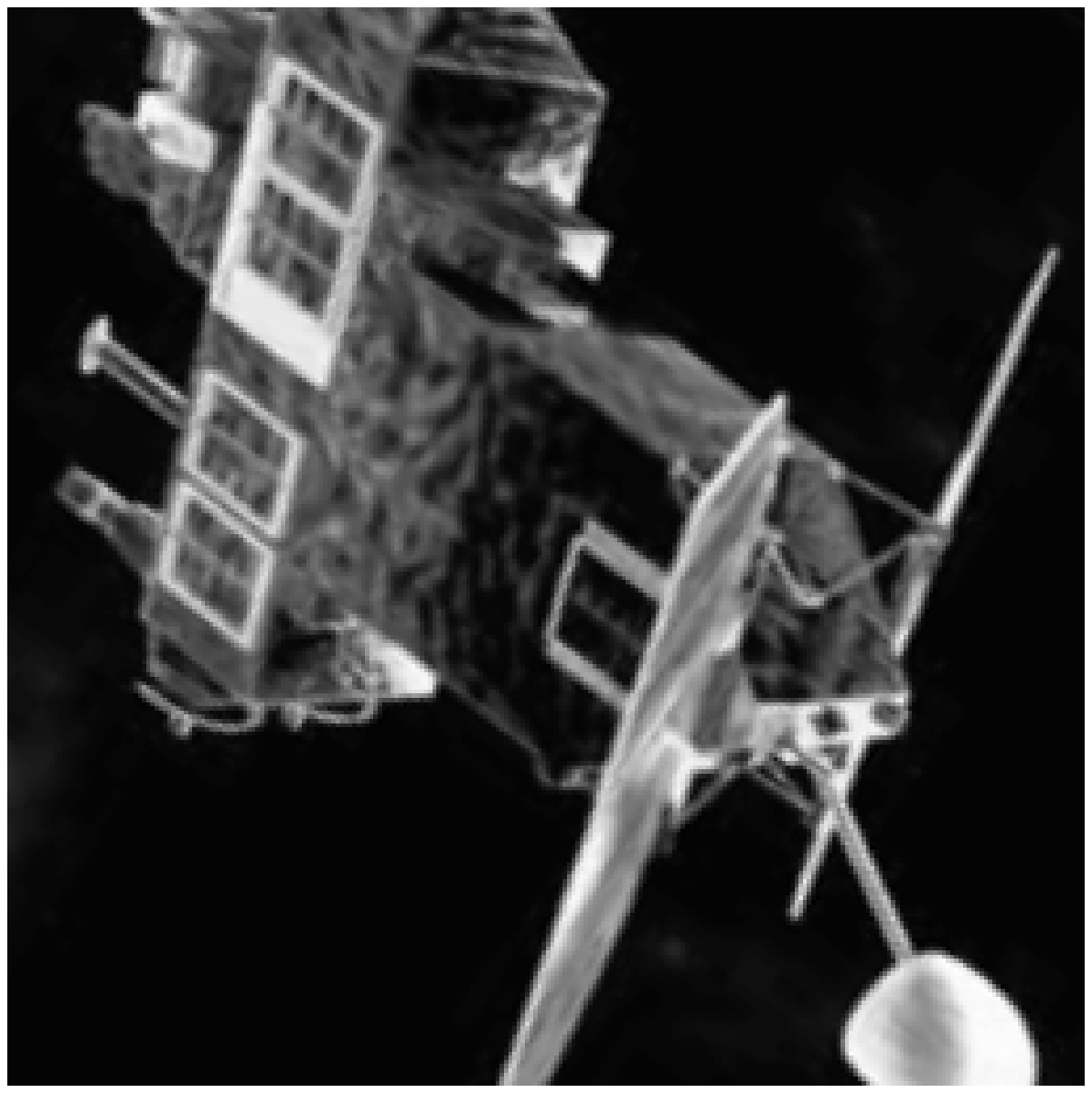}
		\end{tabular}
 	\end{center} 	
 \end{figure}

For the training set, we used eight images of satellites with 10 rigid transformations of each, giving a total of 80 training images.  The blur used was 2D symmetric Gaussian blur with zero mean and variance 1.  Reflexive boundary conditions were assumed. Additive white Gaussian noise was included with a noise level that was randomly selected from a range of .1 and .15.  A set of 80 validation images was generated in the same manner, but with eight different images of satellites.  Sample training images and validation images can be found in the top and bottom rows of Figure~\ref{fig:images} respectively.

We consider various choices for regularization matrix $\bfL$, as well as the multi-parameter Tikhonov case.  We use the identity matrix $\bfL_1 = \bfI$, matrices corresponding to the finite difference approximation to the second derivative in each direction, denoted $\bfL_2$ and  $\bfL_3$, and the discrete Laplacian matrix, $\bfL_4$.  Stencils corresponding to the latter three matrices are given by
\begin{equation}
	\begin{bmatrix} 0 & 0 & 0 \\ 1& -2 & 1\\ 0 & 0 & 0
	\end{bmatrix},\quad
	\begin{bmatrix} 0 & 1 & 0 \\ 0& -2 & 0\\ 0 & 1 & 0
	\end{bmatrix},\quad \mbox{and} \quad
	\begin{bmatrix} 0 & 1 & 0 \\ 1& -4 & 1\\ 0 & 1 & 0
	\end{bmatrix}.
\end{equation}
Since all of these stencils are doubly symmetric and we assumed reflexive boundary conditions, simultaneous diagonalizability~\eqref{eq:simultdiag2} can be obtained where $\bfQ^*$ represents the 2D DCT matrix.

\begin{figure}[bthp]
	\caption{Box plots corresponding to the relative reconstruction errors for the validation set in Example 2.  Results correspond to the 2-norm.}
 	\label{fig:multiLboxplots}	
 	\begin{center}
 		\includegraphics[width=\textwidth]{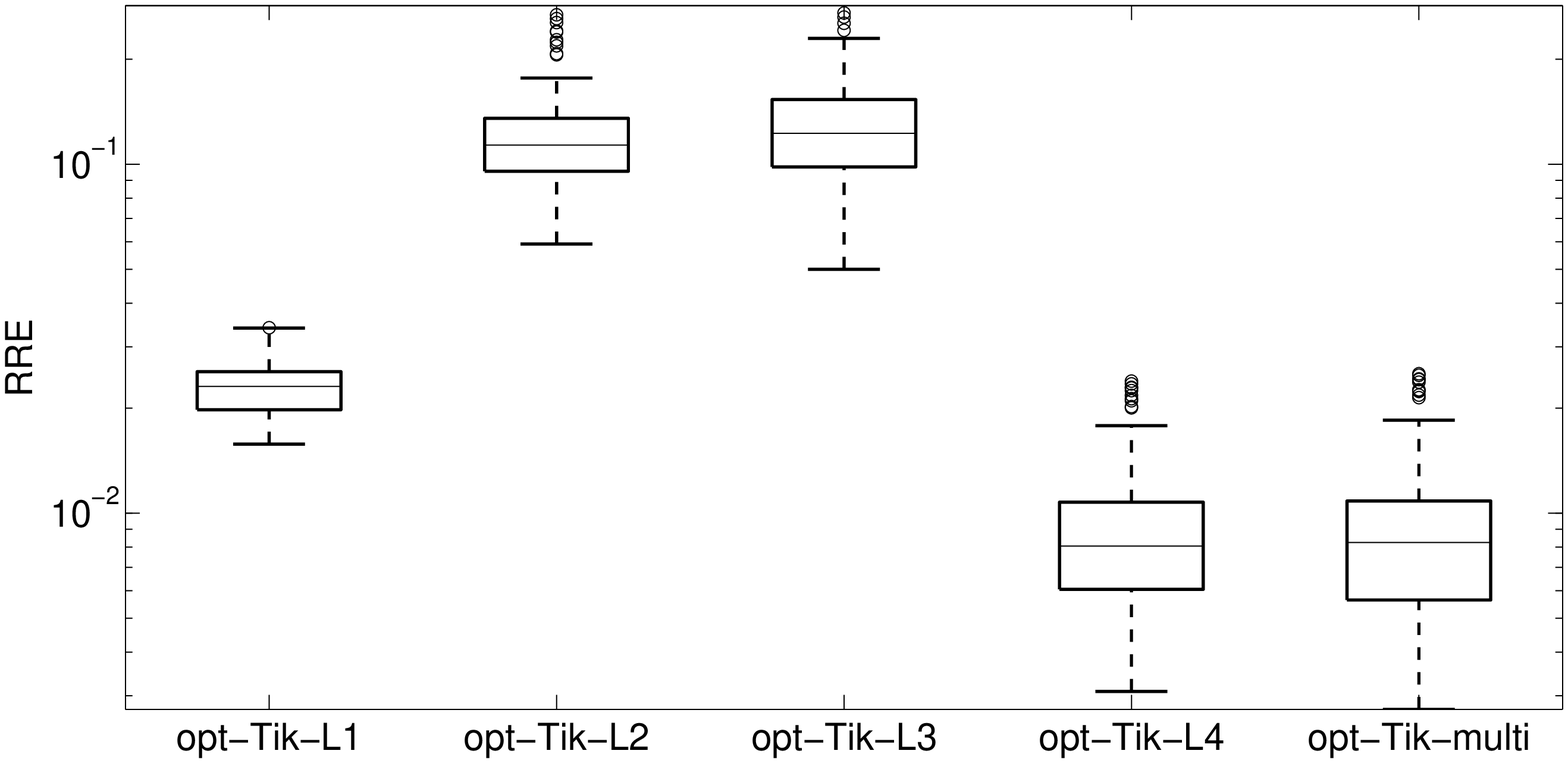}
 	\end{center} 	
 \end{figure}

Using the training data, we compute optimal regularization parameters corresponding to each choice of $\bfL$ individually, which we refer to as opt-Tik-L1, opt-Tik-L2, opt-Tik-L3, and opt-Tik-L4.  In addition, optimal parameters were computed for the multi-parameter Tikhonov problem using all four $\bfL$ matrices; we call this opt-Tik-multi.  Relative reconstruction errors, computed using the 2-norm, for the validation data are presented in Figure~\ref{fig:multiLboxplots} using box plots.   In general, the Laplacian matrix ($\bfL_4$) seems to be a good choice for the regularization matrix, providing smaller relative errors than both standard-form Tikhonov and the one-dimensional derivatives.  Multi-parameter Tikhonov consistently provided good results.

These results correspond to $\rho(\bfxi) = \norm[2]{\bfxi}^2$, but similar results were observed for other error functions.  In particular, we considered $\rho(\bfxi) = \norm[5]{\bfxi}^5$ and $\rho(\bfxi) = \norm[1]{\bfxi}$.  Since optimization of the $1$-norm can be difficult due to non-differentiability, we instead use the Huber function \cite{huber1964robust},
\begin{equation}\rho(\bfxi) =
\begin{cases} |\bfxi| - \frac{\beta}{2}, & \mbox{if } |\bfxi|\geq \beta, \\ \frac{1}{2\beta} \bfxi^2, & \mbox{if } |\bfxi|<\beta, \end{cases}
\end{equation}
with smoothing parameter $\beta = 10^{-4}$, which provides an approximation to the $1$-norm.
Relative reconstruction errors were computed as $\rho(\bfxi_\lambda - \bfx_\true)/\rho( \bfx_\true)$. The average relative reconstruction error for the validation set and corresponding standard deviations, can be found in Table~\ref{tab:MultiLerrors} for the different approaches.  In addition, we compare to GCV-Tik-multi, which corresponds to multi-parameter Tikhonov, where the regularization parameter vector $\bflambda$ minimizes the GCV function in~\eqref{eq:GCVmultiGSVD}.  The parameter vector is problem dependent and was computed separately for each validation image.

\begin{table}[bthp]
	\caption{Average relative reconstruction error and standard deviation, in parenthesis, for the validation set in Example 2 for various error functions.}
	 \label{tab:MultiLerrors}	 		
	 	\begin{center}
    \scalebox{0.9}{
	\begin{tabular}{|l||c|c|c|} \hline
   	&  Huber   & 2-norm & 5-norm \\ \hline
opt-Tik-$L_1$  & 1.672e-01 (4.164e-02) & 2.324e-02 (4.375e-03) &  1.232e-04 (9.215e-05) \\ \hline
opt-Tik-$L_2$  & 3.718e-01 (1.411e-01) & 1.257e-01 (4.904e-02) &  1.118e-02 (9.886e-03) \\ \hline
opt-Tik-$L_3$  & 3.816e-01 (1.422e-01)  & 1.303e-01 (4.669e-02)    & 1.145e-02 (9.391e-03)  \\  \hline
opt-Tik-$L_4$ & 9.368e-02 (4.586e-02)  & 9.766e-03 (5.457e-03)   & 3.383e-05 (3.436e-05)  \\  \hline
opt-Tik-multi  & 9.393e-02 (4.806e-02)  & 9.829e-03 (5.993e-03)  & 3.469e-05 (3.817e-05)  \\  \hline
GCV-Tik-multi  & 9.411e-02 (4.739e-02) & 8.870e-03 (5.172e-03)  & 5.076e-05 (4.711e-05)  \\  \hline
\end{tabular} }
\end{center}

\end{table}

We observe that in all error measures, the Laplacian matrix outperforms all of the tested $\bfL$ matrices, even in this case, the multi-parameter Tikhonov case.  Note that in the 2-norm, GCV-Tik-multi performs very well, even better than opt-Tik-multi.  This is expected since opt-Tik-multi computes the optimal regularization parameters \emph{on average}, while GCV-Tik-multi computes regularization parameters that are specific to each problem. However, for the other error measures, this is not the case and opt-Tik-multi results in smaller average reconstruction errors than GCV-Tik-multi.

\begin{figure}[bthp]
	 \caption{Absolute error images, in inverted colormap so that white corresponds to small errors, for one of the validation images in Example 2. Images on the left correspond to optimal parameters that were computed using training data. The two images on the right correspond to GCV-selected parameters for both standard Tikhonov and multi-parameter Tikhonov.}
 	\label{fig:multiLerrorimg}
 	\begin{center}
		\begin{tabular}{cccc|cc} & Huber & 2-norm & 5-norm & GCV\\
\rotatebox{90}{opt-Tik-$\bfL_1$}
	&	\includegraphics[scale=0.2]{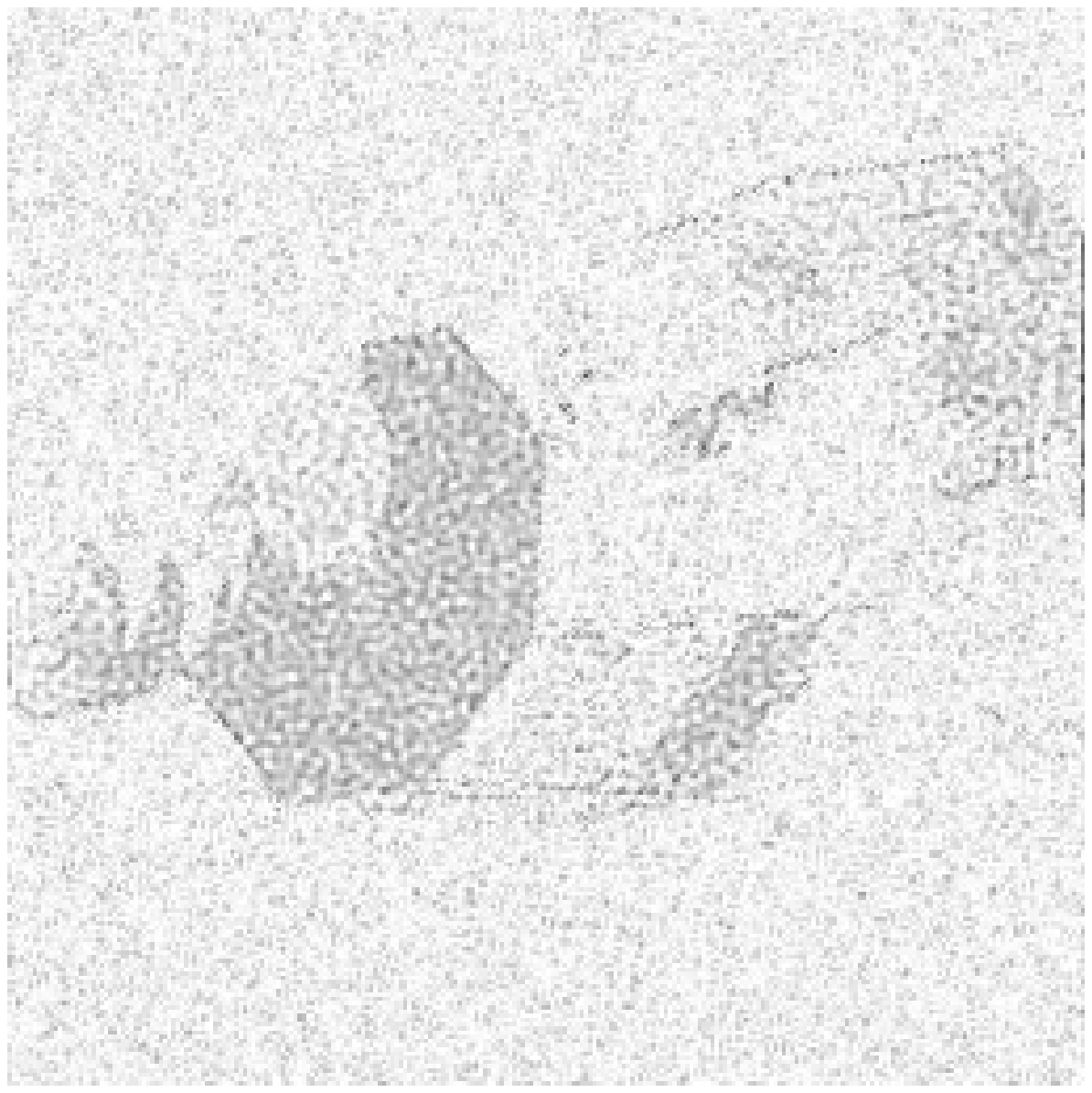} &
 		\includegraphics[scale=0.2]{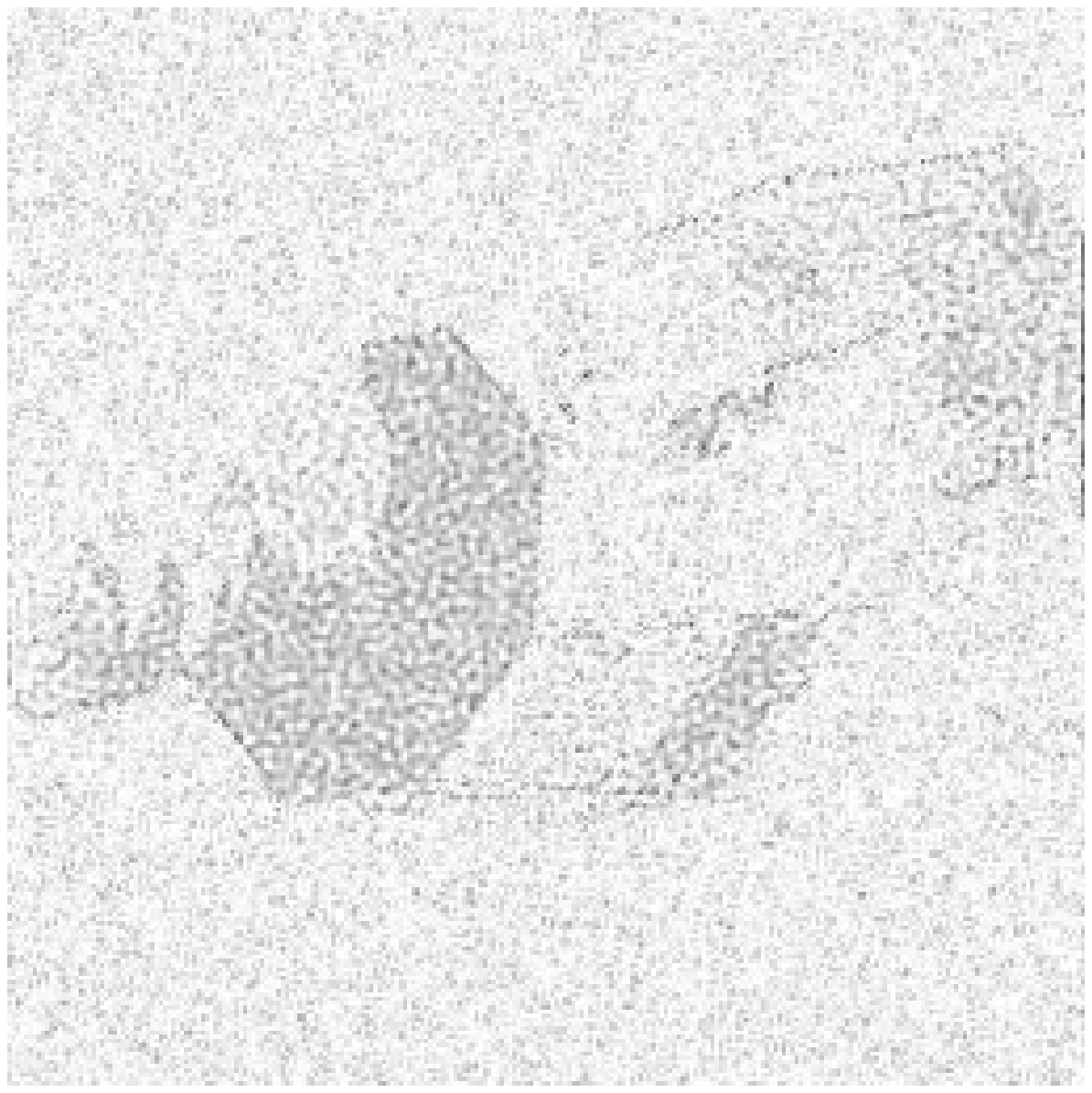} &
 		\includegraphics[scale=0.2]{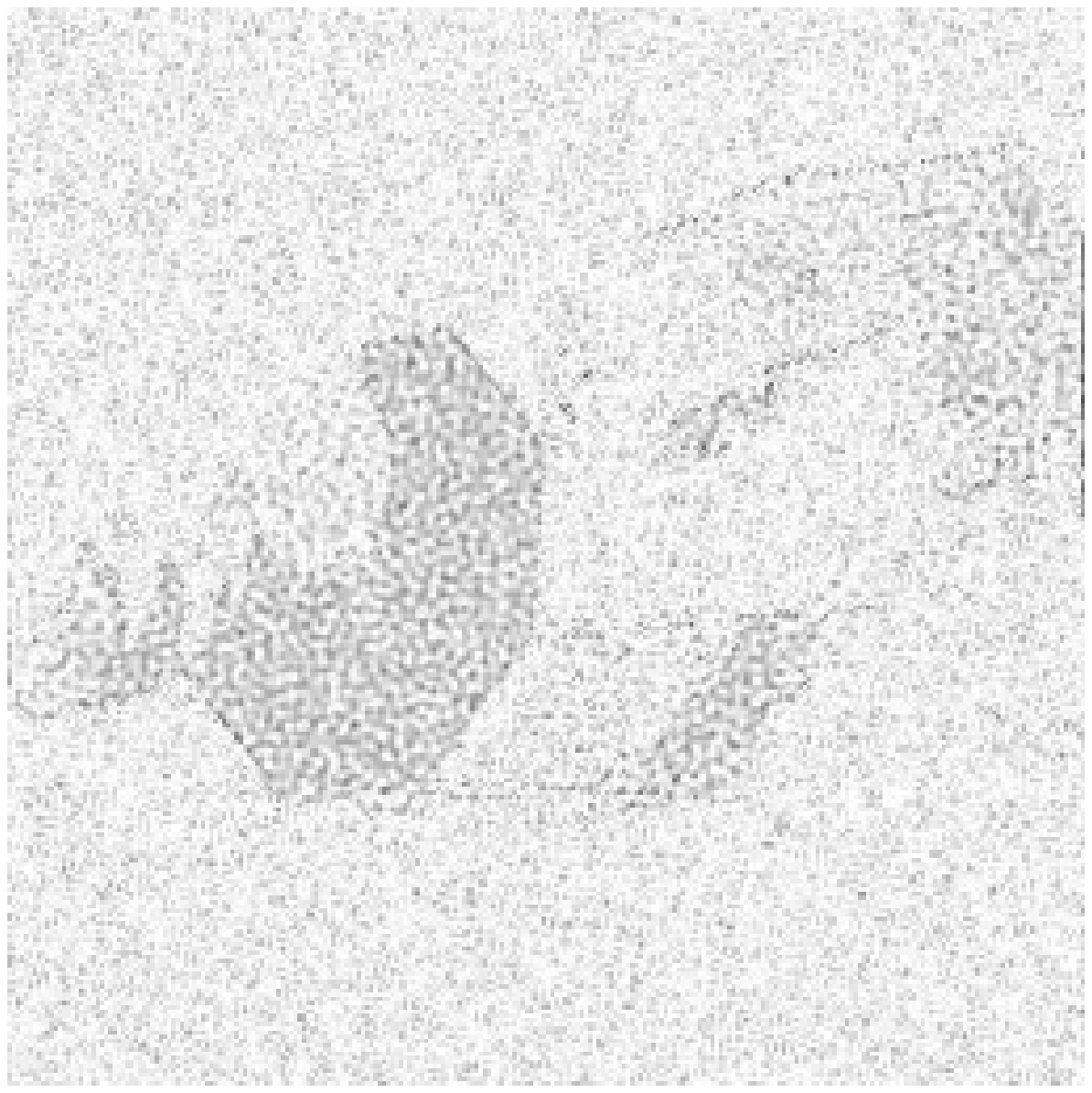} &
 		\includegraphics[scale=0.2]{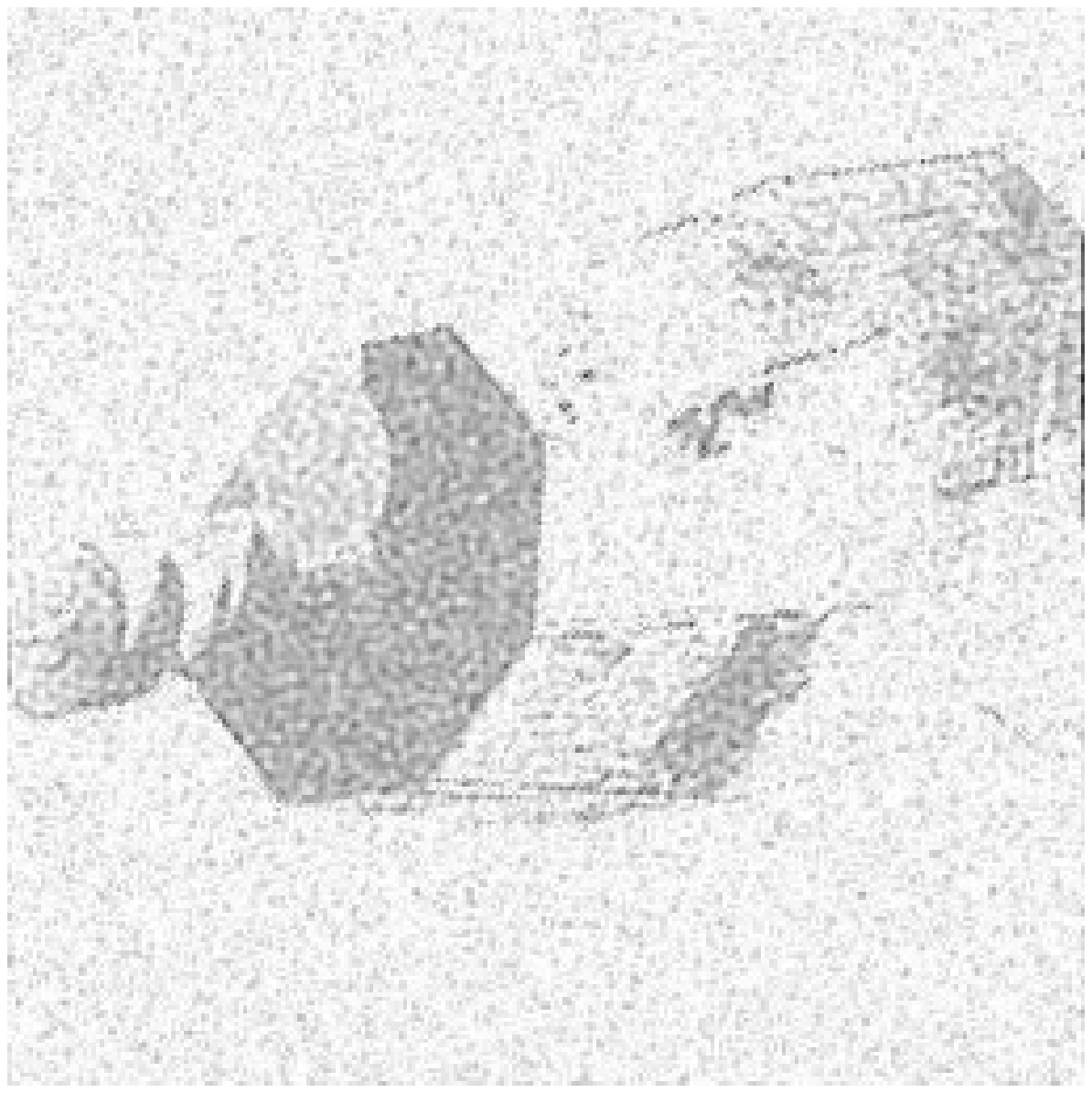} & \rotatebox{90}{Tik-$\bfL_1$}\\
 \rotatebox{90}{opt-Tik-multi}	 &	\includegraphics[scale=0.2]{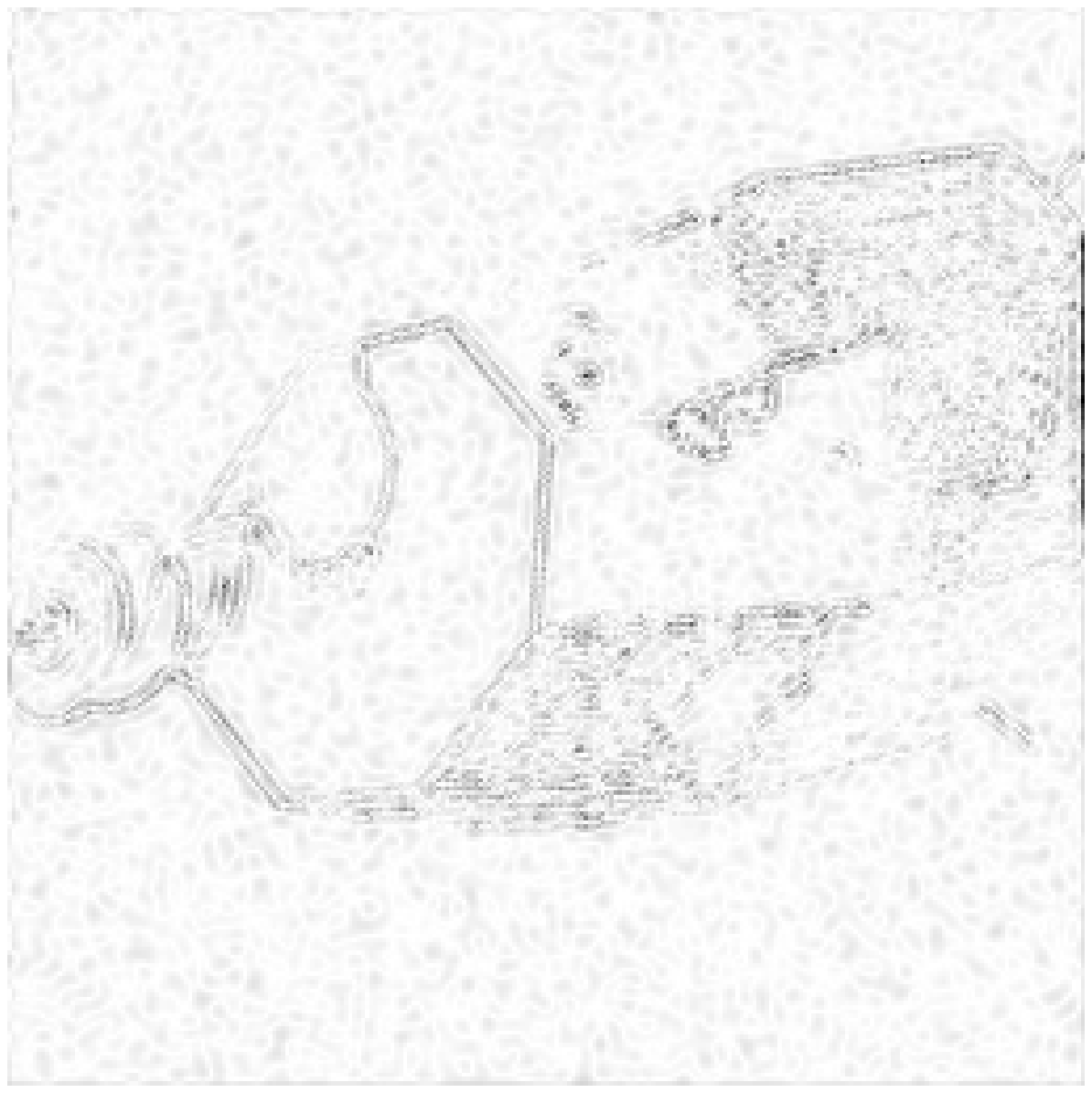} &
 		\includegraphics[scale=0.2]{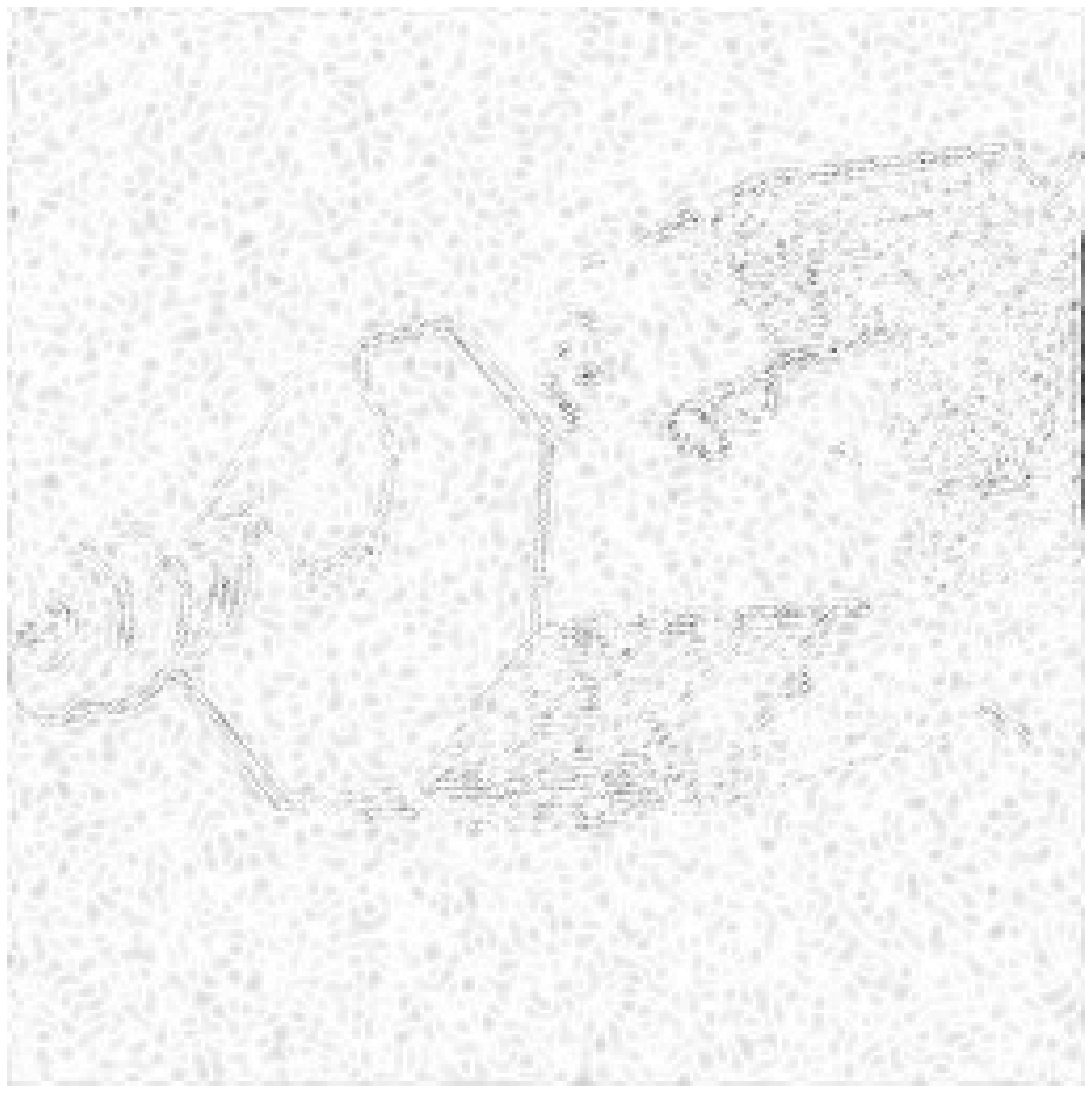} &
 		\includegraphics[scale=0.2]{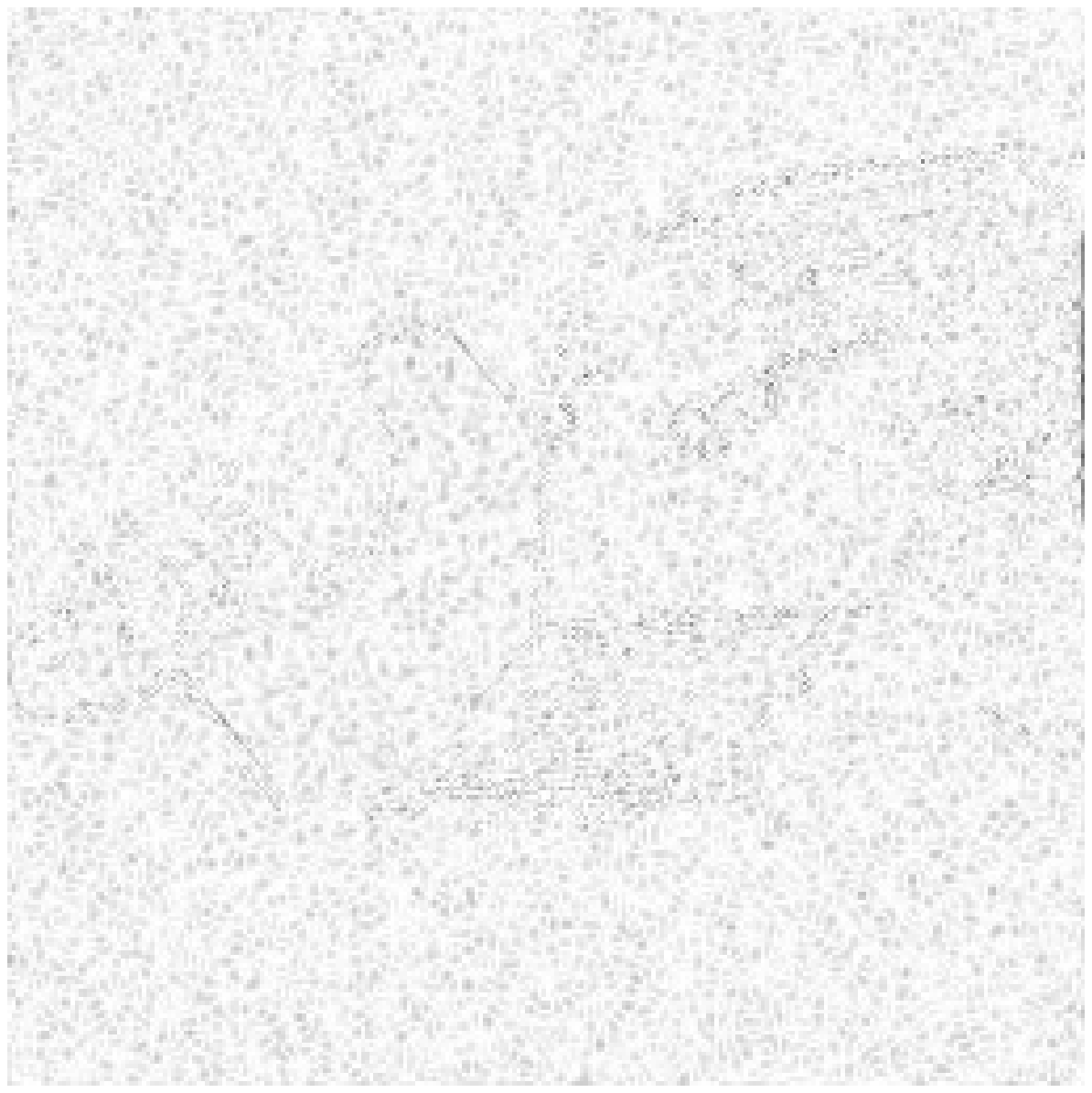} &
 		\includegraphics[scale=0.2]{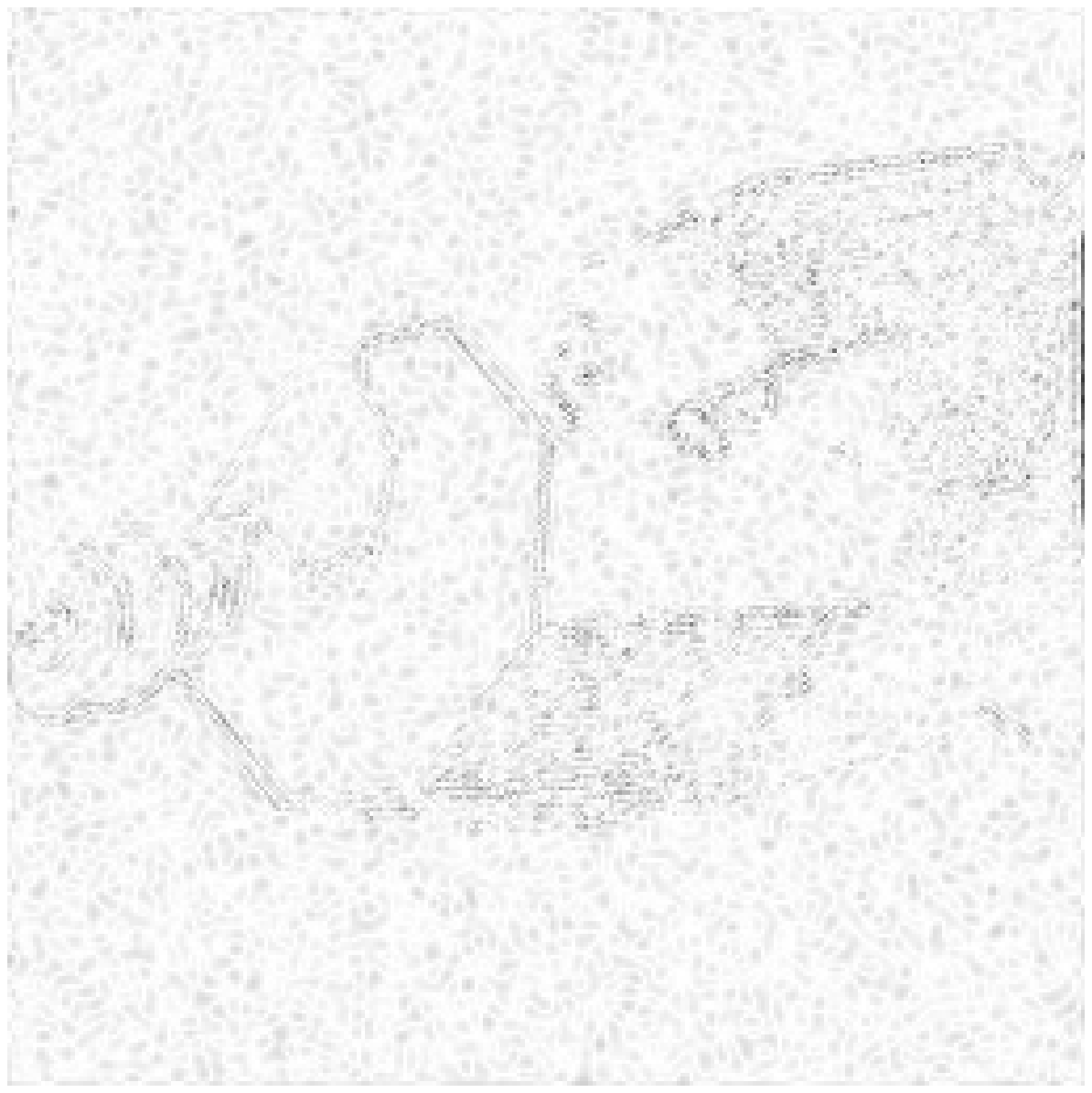} &  \rotatebox{90}{Tik-multi}
		\end{tabular}
 	\end{center}	
 \end{figure}

The specific choice of $\rho$ is reflected in the quality of the reconstructed images.  For one validation image, we present in Figure~\ref{fig:multiLerrorimg} absolute error images (i.e., absolute value of the reconstruction minus the true image) in inverted colormap for opt-Tik-$L_1$ (standard Tikhonov) and opt-Tik-multi.  Similar to the above conclusions, we see that opt-Tik-multi provides reconstructions with smaller absolute errors than standard Tikhonov.  Comparing the different error functions, we observe that reconstructions with the Huber function have small errors in flat or constant regions of the image and large errors near edges.  Furthermore, the 5-norm reconstruction has overall smaller reconstruction errors, but the errors are more distributed in the image. This is consistent with observations made in \cite{chung2011designing}.  The error images on the far right correspond to GCV for both standard Tikhonov and multi-parameter Tikhonov.

 \paragraph{Example 3}
In this example we consider a multi-parameter Tikhonov problem where operator approximations are used to estimate regularization parameters.  As mentioned in Section~\ref{sub:multiple_parameter_tikhonov}, for problems where simultaneous diagonalizability does not exist or is not easily obtainable, for example for spatially variant blurs \cite{nagy1998restoring} or problems where $\mxL$ represents dictionaries or gradient masks \cite{liu2012morphology}, a bilevel optimization approach could be used to compute optimal regularization parameters, but such an approach may be costly.  We follow the approach described in \cite{ChungKilmerOleary} where operator approximations are used to estimate regularization parameters, and then the original problem is solved using these parameters.

For this example, we use the 2D image deblurring problem described in Example 2, where $\bfA$, $\bfL_1$, $\bfL_2$, $\bfL_3$ and $\bfL_4$ are jointly diagonalizable by the 2D DCT due to the assumption of reflexive boundary conditions.  However, we assume for the moment that such a factorization is not available, and we use the assumption of periodic boundary conditions to obtain approximate matrices $\bf{\hat{A}}$, $\hat \bfL_1$, $\hat \bfL_2$, $\hat \bfL_3$ and $\hat \bfL_4$, which are jointly diagonalizable by the 2D DFT. Using the training data, we compute the optimal regularization parameters corresponding to the approximate problem, denoted as $\hat \bflambda$, by using the Gauss-Newton approach described in Section \ref{sub:multiple_parameter_tikhonov}.  For comparison purposes,  we also use the training data to compute the optimal regularization parameters for the original problem.

\begin{figure}
         \caption{Box plots of relative reconstruction errors for opt-Tik-multi, Tik-multi-$\hat\bflambda$, and Tik-multi-$\hat \bflambda_\GCV$.}
 	\label{fig:example3_box}
 \centering
 \includegraphics[width=\textwidth]{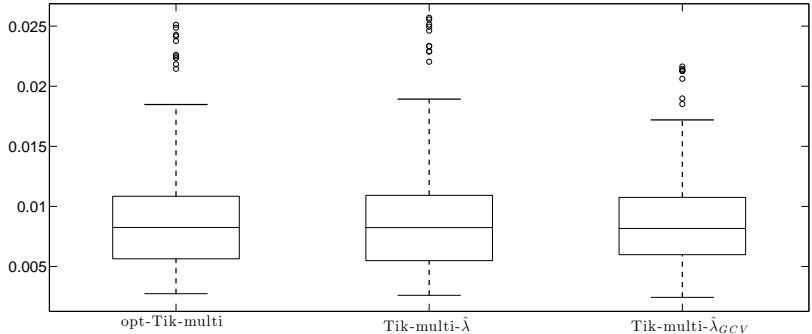}
 \end{figure}

Regularization parameters were then used to reconstruct the blurred images in the validation set.  Results using optimal regularization parameters corresponding to the original problem are provided as opt-Tik-multi (this corresponds to Figure~\ref{fig:multiLboxplots}), while Tik-multi-$\hat\bflambda$ corresponds to reconstructions from solving the original problem with regularization parameters $\hat \bflambda$.  For each validation image, we also use the multi-parameter GCV function on the approximate problem to estimate parameters $\hat \bflambda_\GCV$ that were then used to solve the original problem.  This approach is denoted Tik-multi-$\hat \bflambda_\GCV$.
We compute the relative 2-norms of the errors for the reconstructions of the validation data and display their distributions using boxplots in Figure~\ref{fig:example3_box}.

From these comparisons, we see optimal parameters computed for the approximate problem can lead to relative reconstruction errors that are comparable to errors obtained using optimal parameters for the original problem.  An interesting phenomenon is that the relative reconstruction errors for Tik-multi-$\hat \bflambda_\GCV$ have a slightly smaller interquartile range and smaller outliers. This is because GCV parameters are computed individually for each image in the validation set thereby tailoring the parameters to the individual characteristics of each image in the validation set.  On the other hand, the multi-parameter approach takes a holistic approach and finds the best set of regularization parameters that minimizes errors for the entire set.

\begin{figure}
        \caption{Histograms of GCV-computed regularization parameters $\hat \bflambda_\GCV$ over all validation images.  Solid vertical lines corresponds to computed parameters $\hat \lambda$ that were obtained using approximate operators and training data.}
	\label{fig:lambdaDist}
\centering
\begin{tabular}{cc}
\includegraphics[width=.4\textwidth]{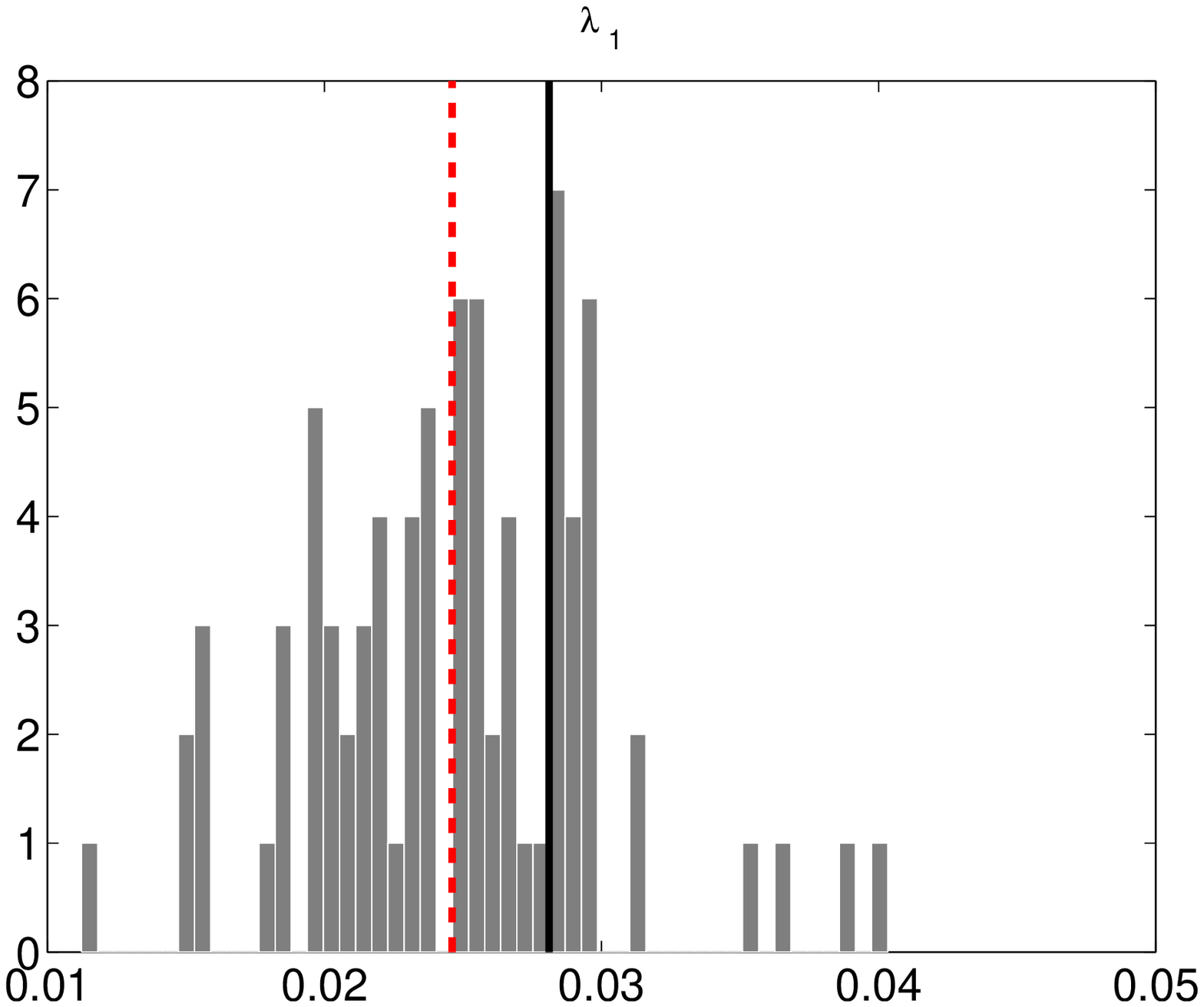} &
\includegraphics[width=.4\textwidth]{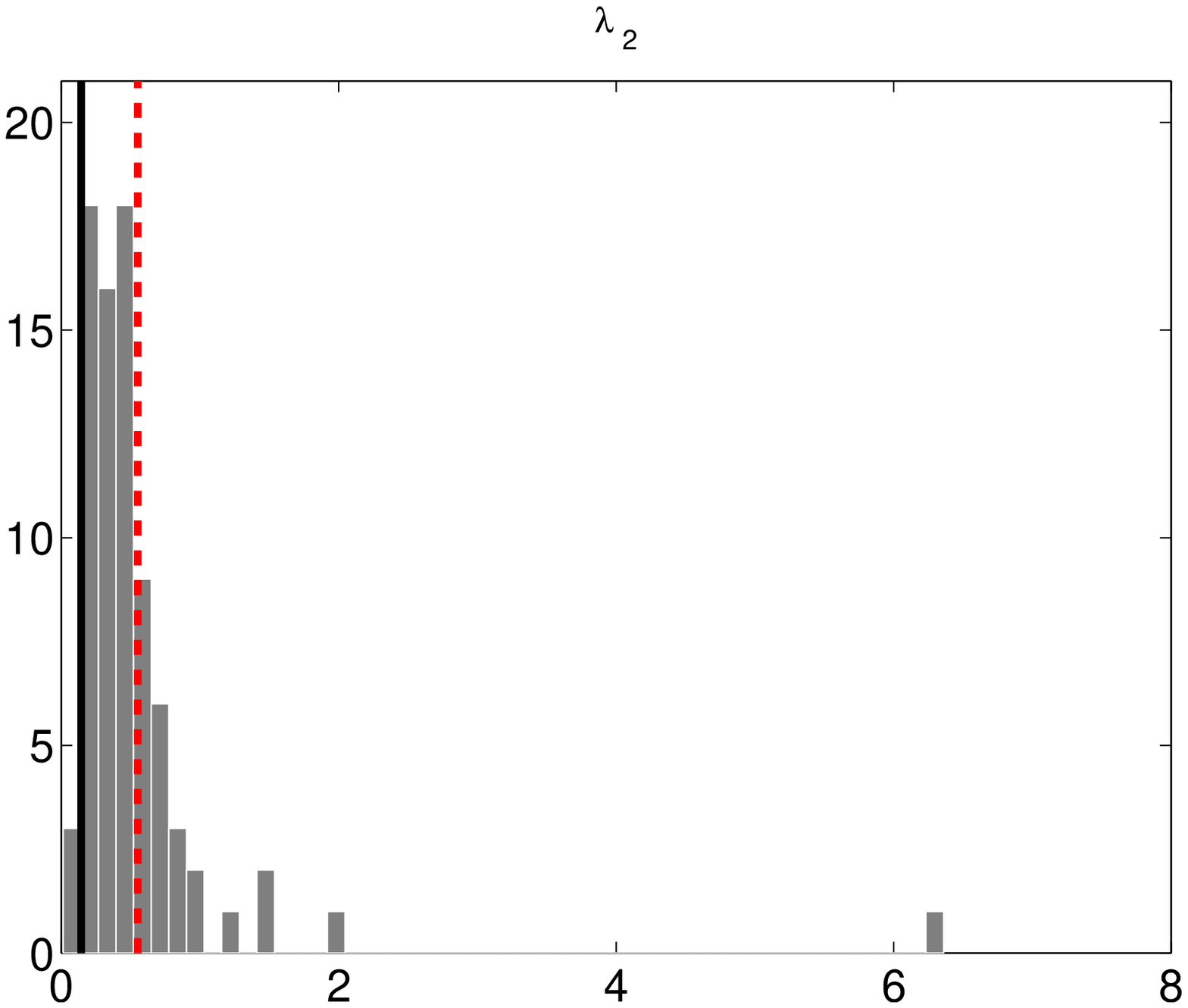} \\
\includegraphics[width=.4\textwidth]{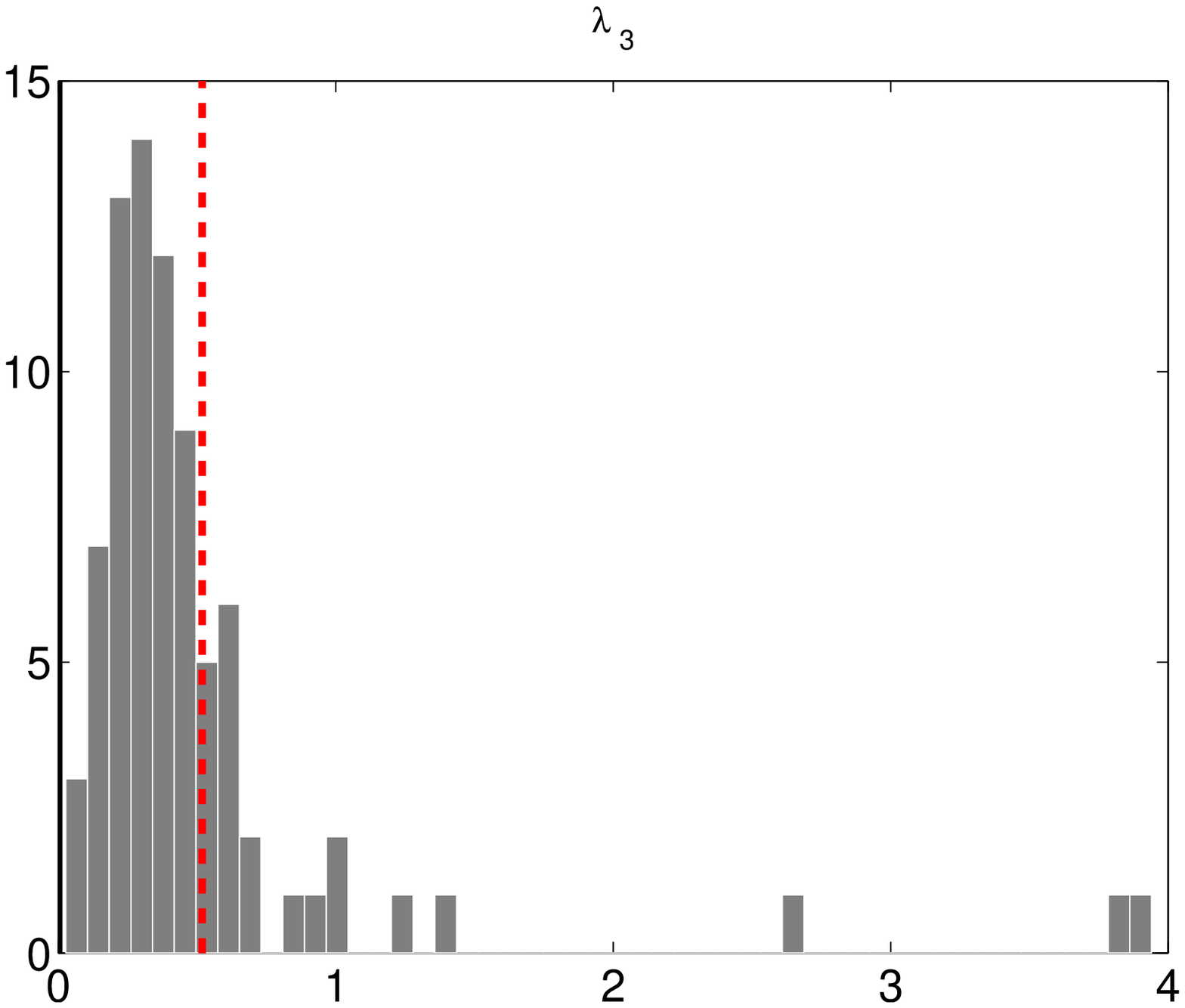} &
\includegraphics[width=.4\textwidth]{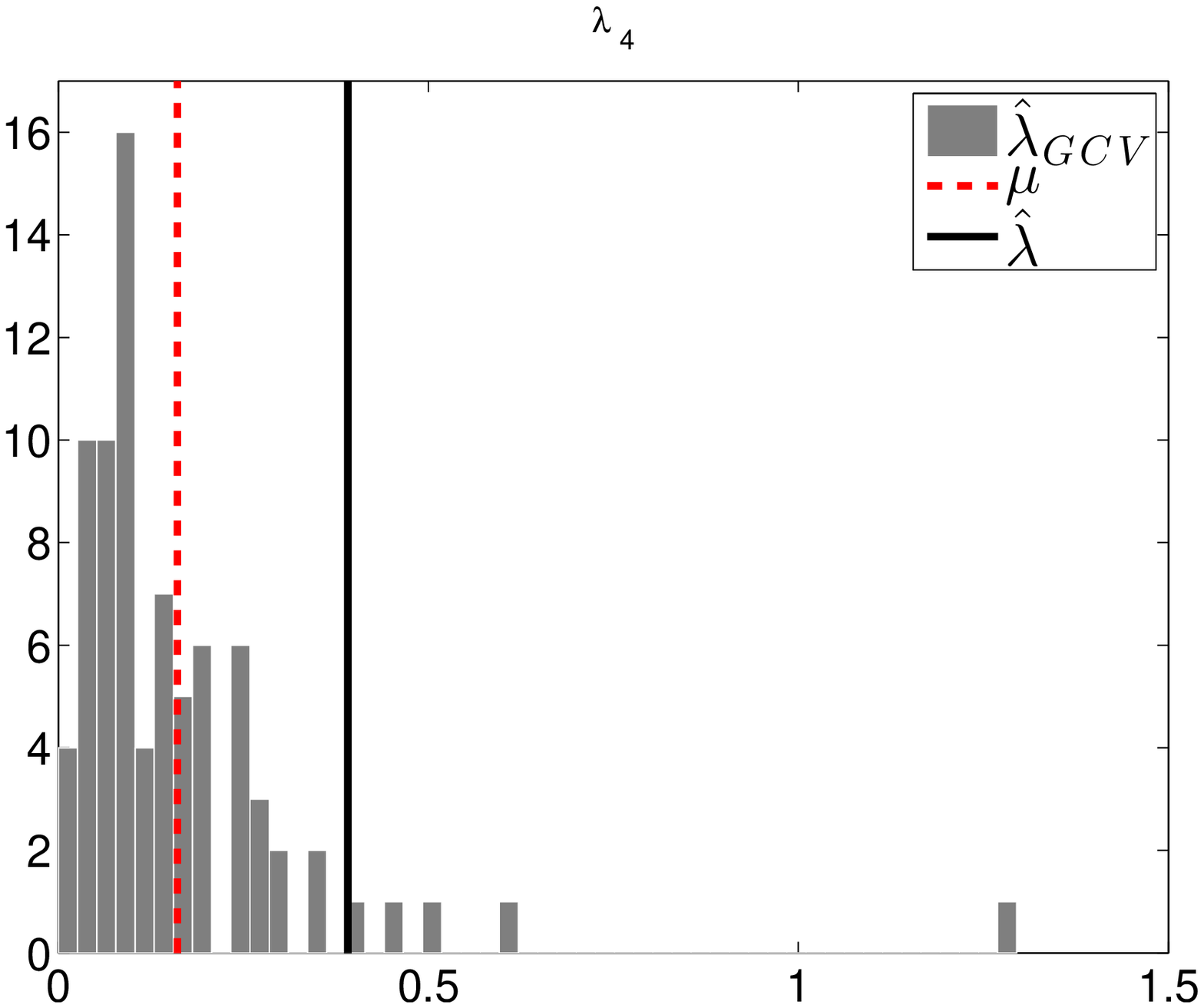}
\end{tabular}
\end{figure}

We investigate this further by considering the values of the computed regularization parameters $\hat \bflambda$ and $\hat \bflambda_\GCV$.  In Figure~\ref{fig:lambdaDist},  we provide histograms of each of the GCV-computed parameters over all validation images, with the corresponding mean value $\mu$.  In addition, we provide the computed parameters $\hat \bflambda$ as solid vertical lines.
As seen in the histogram, Tik-multi-$\hat \bflambda$ overweights the significance of $\bf{L}_1$ and $\bf{L}_4$ and underweights that of $\bf{L}_2$ and  $\bf{L}_3$. This can be explained by the fact that $\bf{L}_1$ and $\bf{L}_4$ take into account information in both vertical and horizontal directions, while $\bf{L}_2$ and $\bf{L}_3$ only look in one direction, which may be beneficial for only a few images. It is worthwhile to mention the computational trade-off between the two approaches.  In general, using GCV on each individual image can produce superior results, but can be costly for problems where very large sets of images must be reconstructed in real time.  For these problems, finding one set of parameters that will work for a variety of images (e.g., opt-Tik-multi or multi-Tik-$\hat\bflambda$) will have a significant cost advantage and will produce results that are nearly as accurate as those obtained by individual GCV results.

In these examples, satellite images (although different ones) were used for both training and validation. A natural question would be to investigate how different the images can be from the training set and still produce good reconstructions. That is, holding everything else constant, can the regularization parameters that were computed for a set of satellite images be used to reconstruct other types of images?

\begin{figure}
	 \caption{True, blurred, and reconstructed images for the Cell, Fabric, Ariel Map, and Camera Man examples.  Reconstructions were computed using multi-parameter Tikhonov regularization, where the regularization parameters correspond to the training set of satellite images. }
 	\label{fig:example3_Outside}
 	\begin{center}
		\begin{tabular}{ccc}
         True & Blurred & Reconstructed\\
	\includegraphics[scale=0.2]{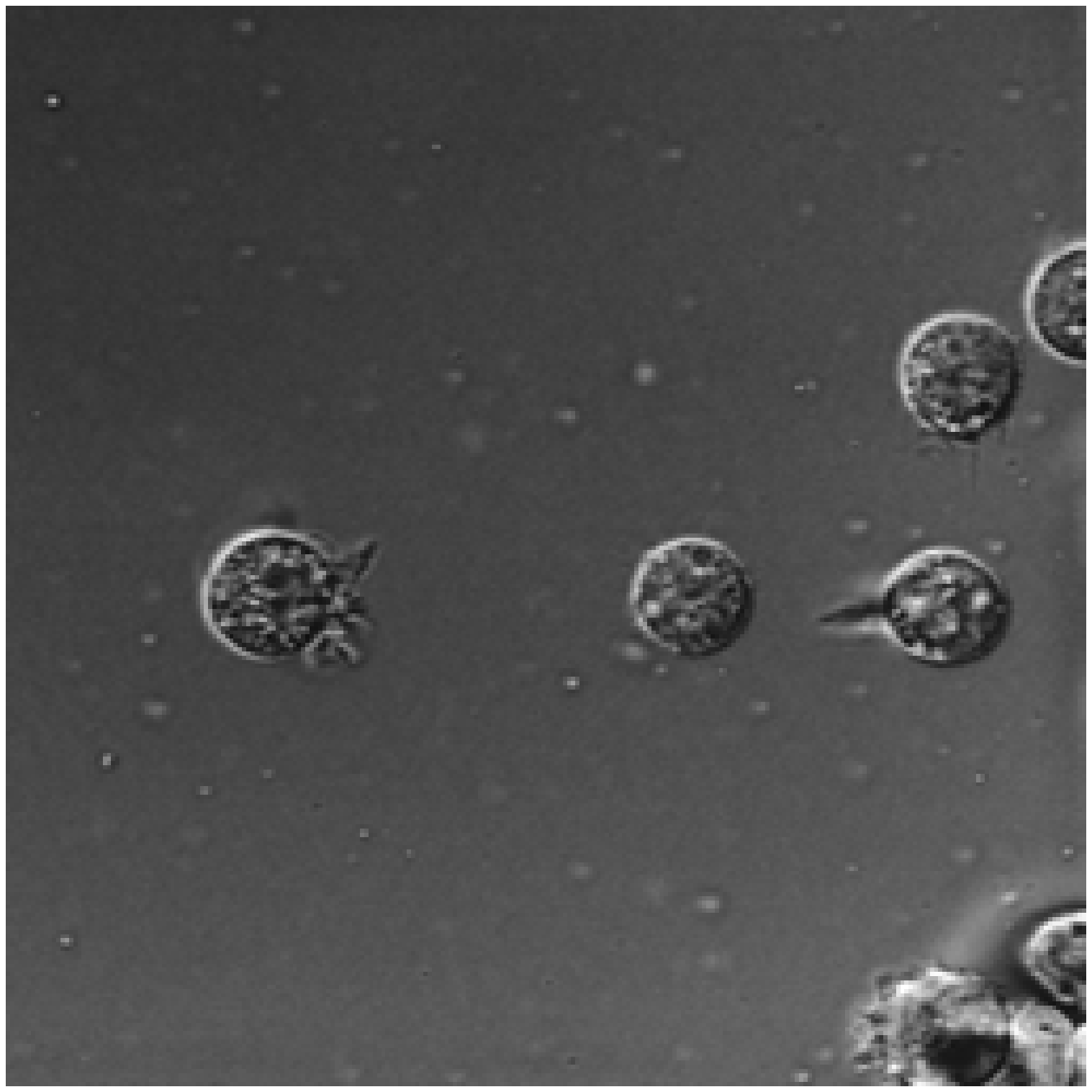} &
 		\includegraphics[scale=0.2]{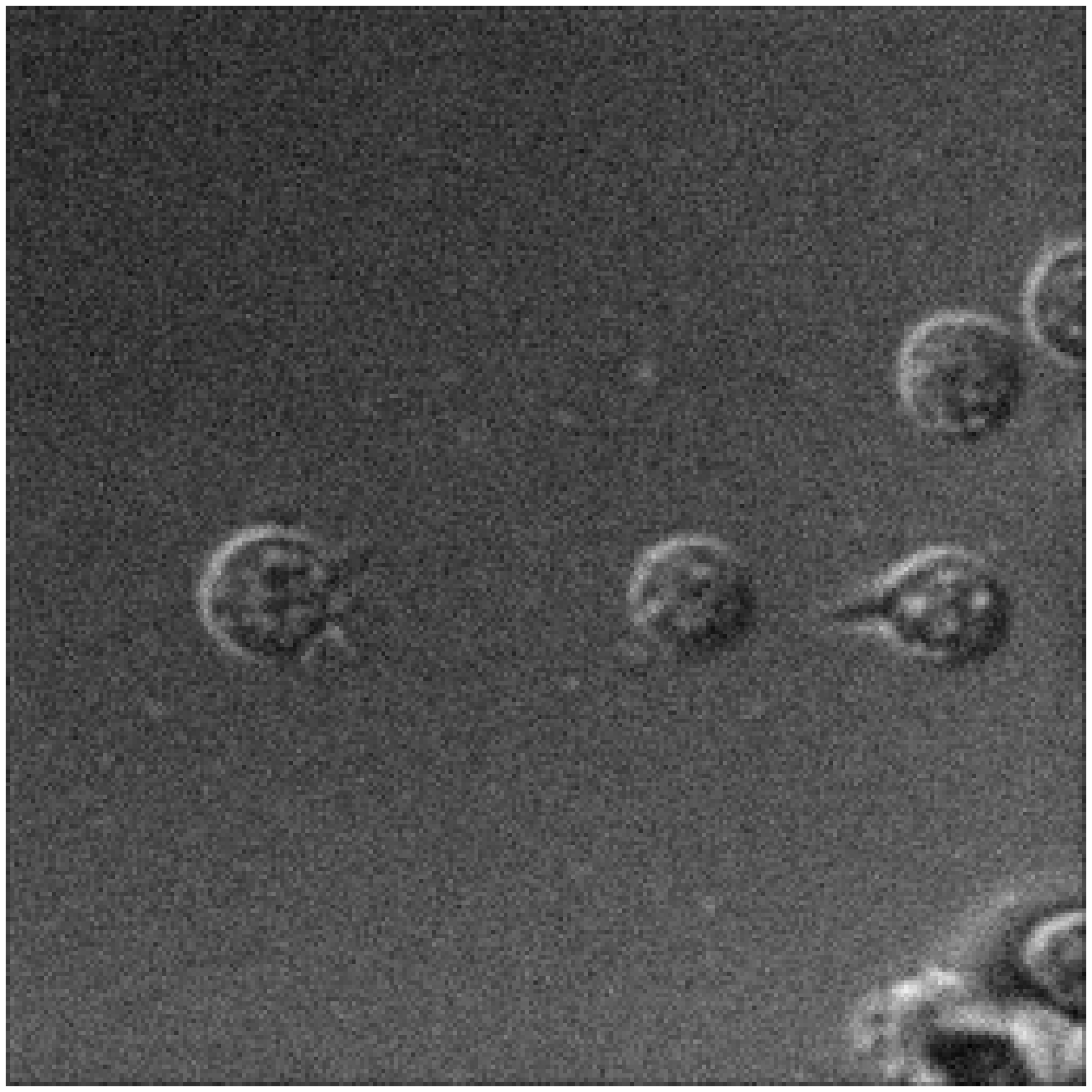} &
 		\includegraphics[scale=0.2]{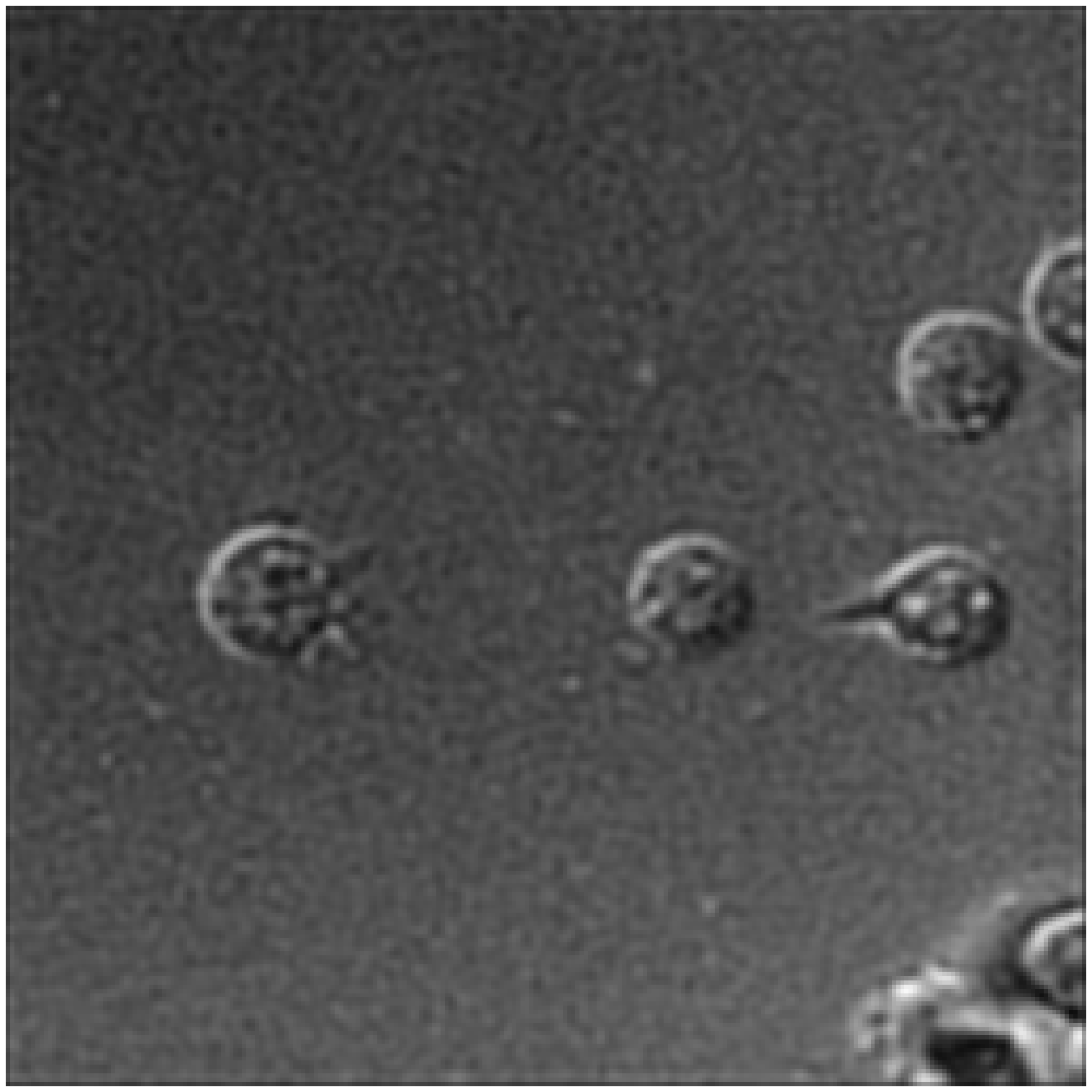} \\
\includegraphics[scale=0.2]{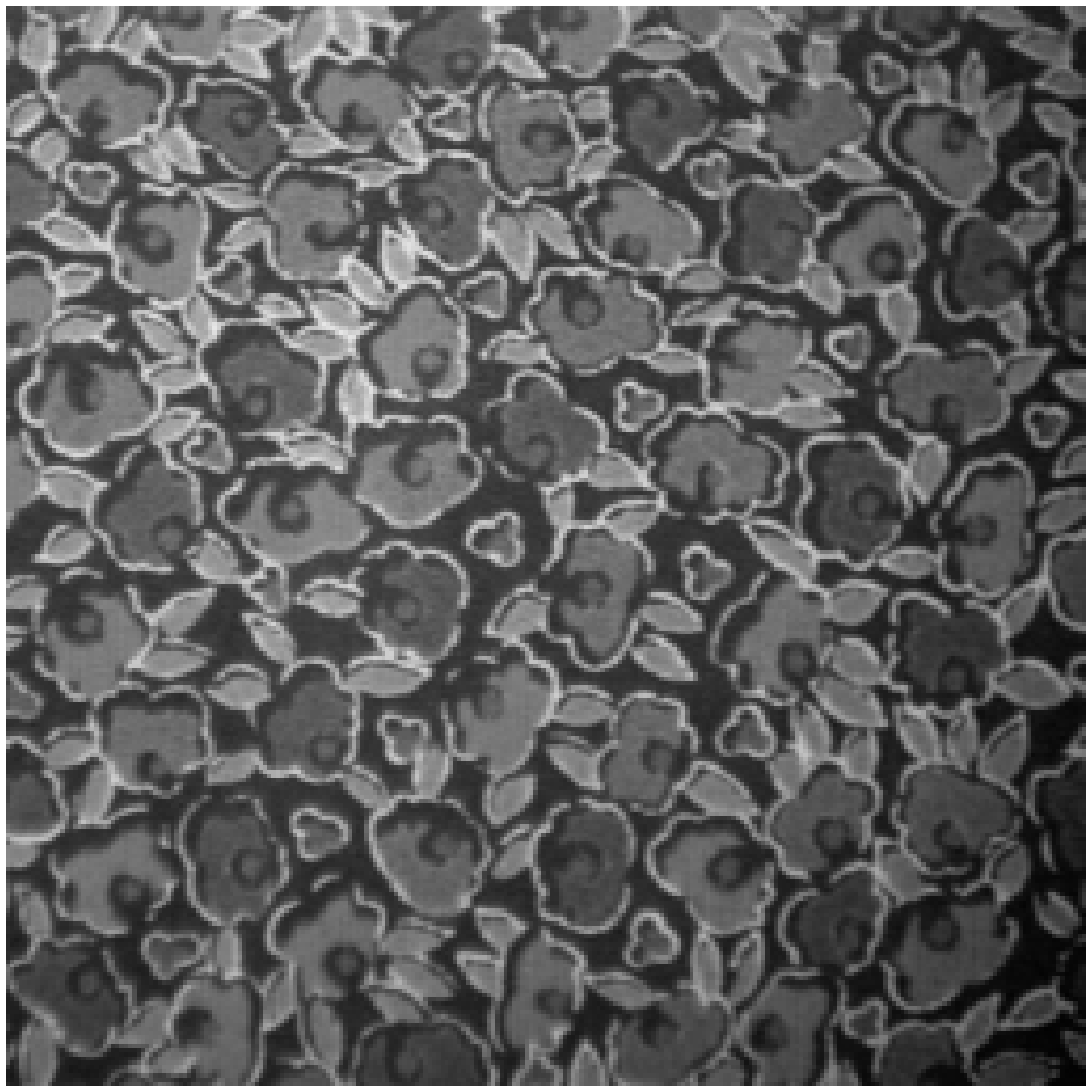} &
 		\includegraphics[scale=0.2]{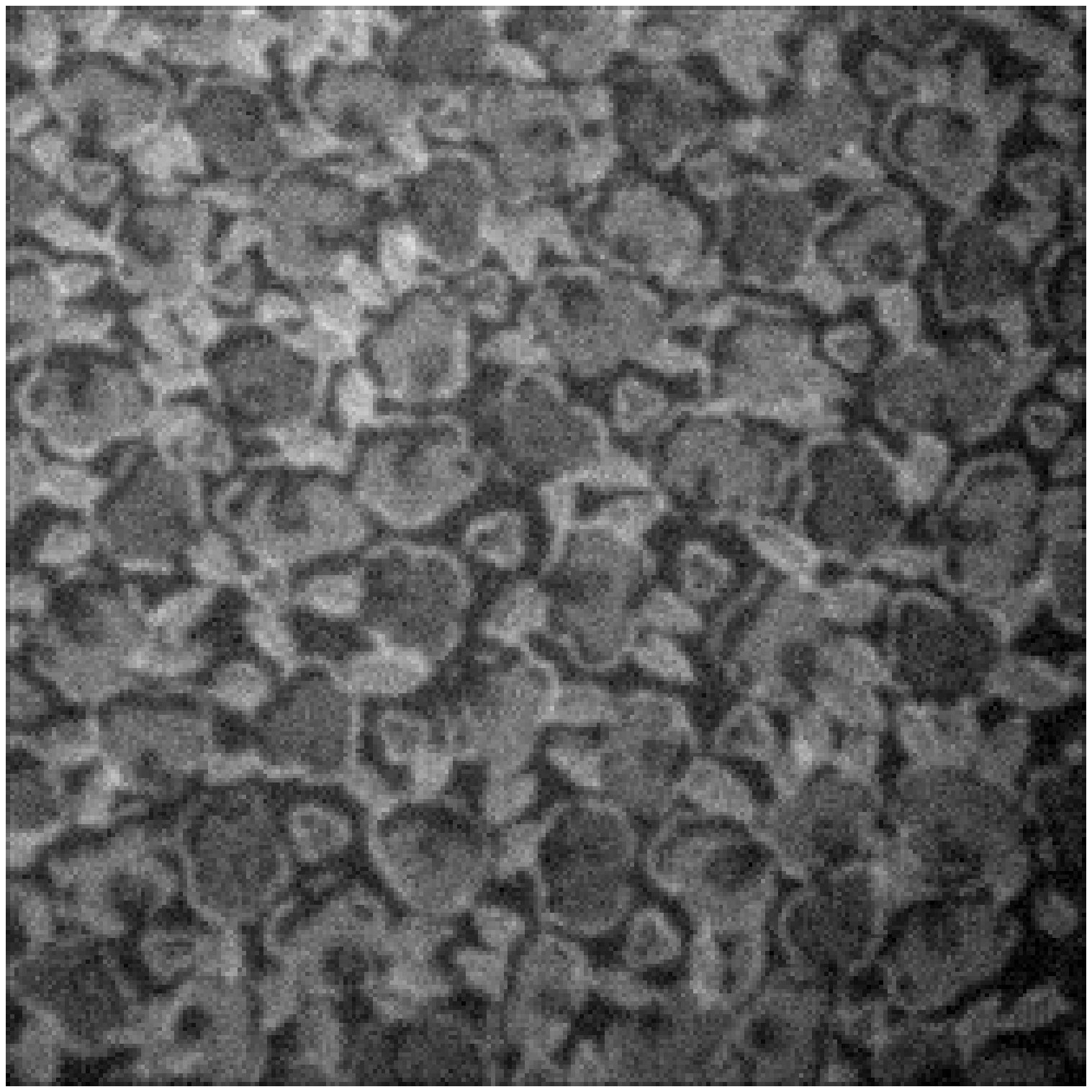} &
 		\includegraphics[scale=0.2]{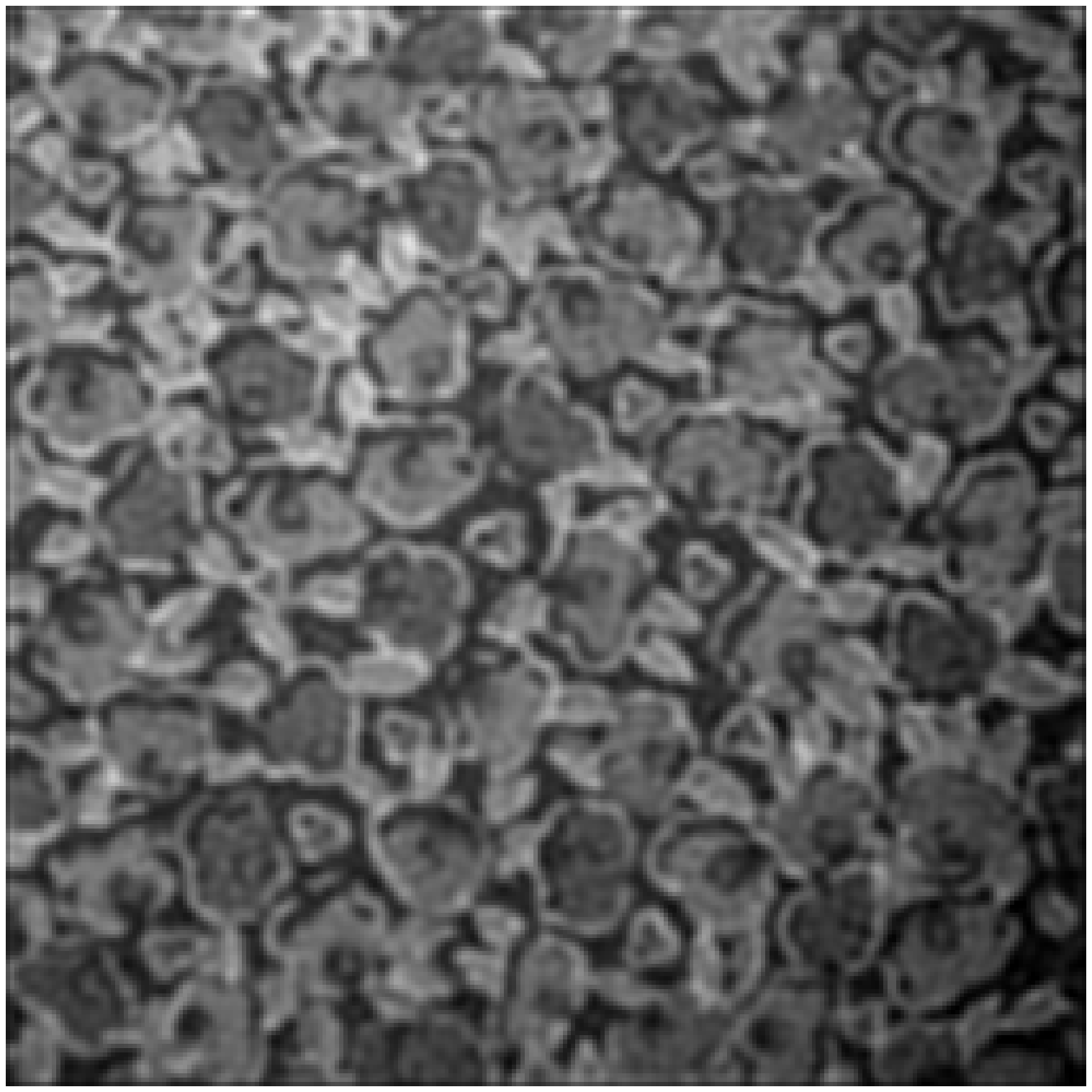} \\
\includegraphics[scale=0.2]{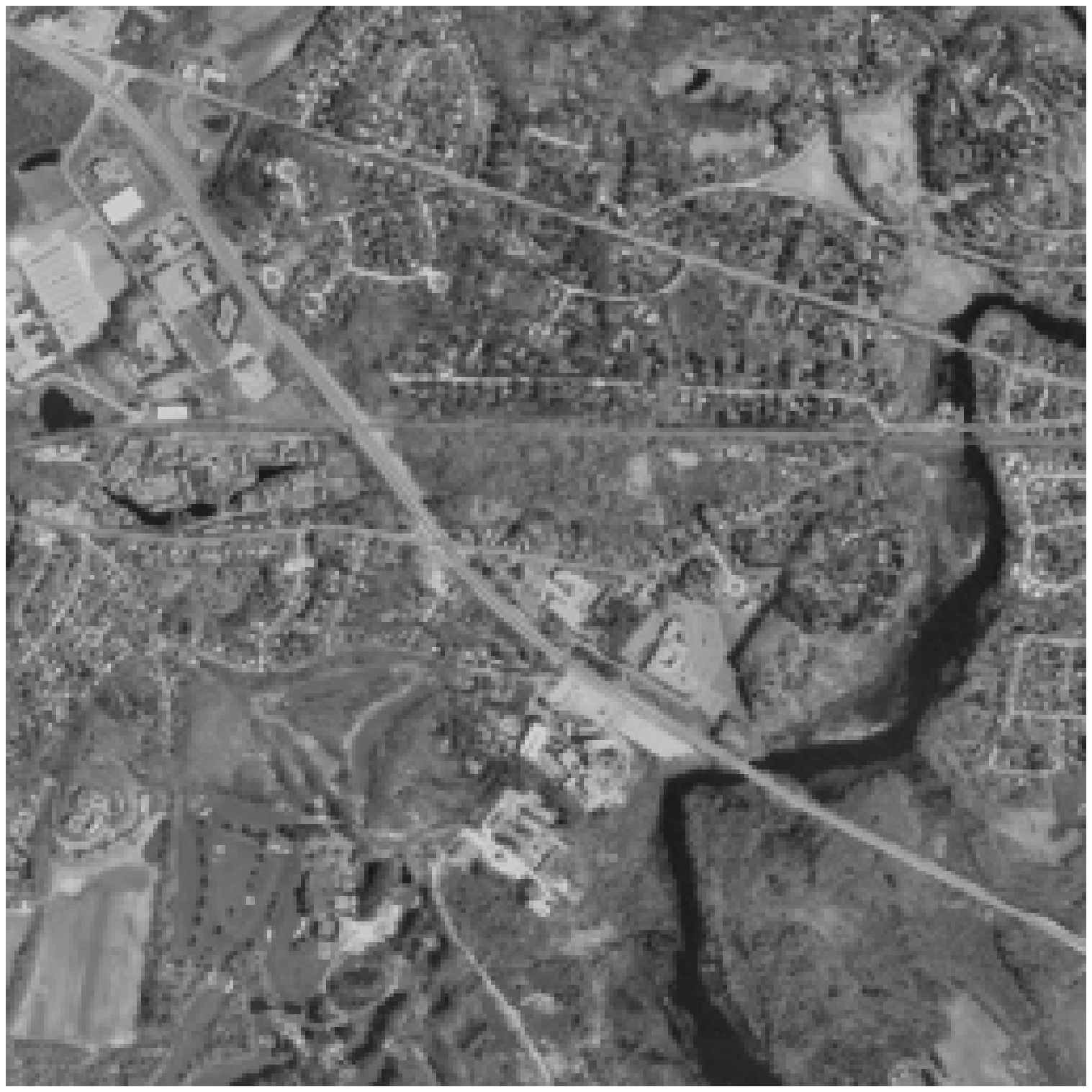} &
 		\includegraphics[scale=0.2]{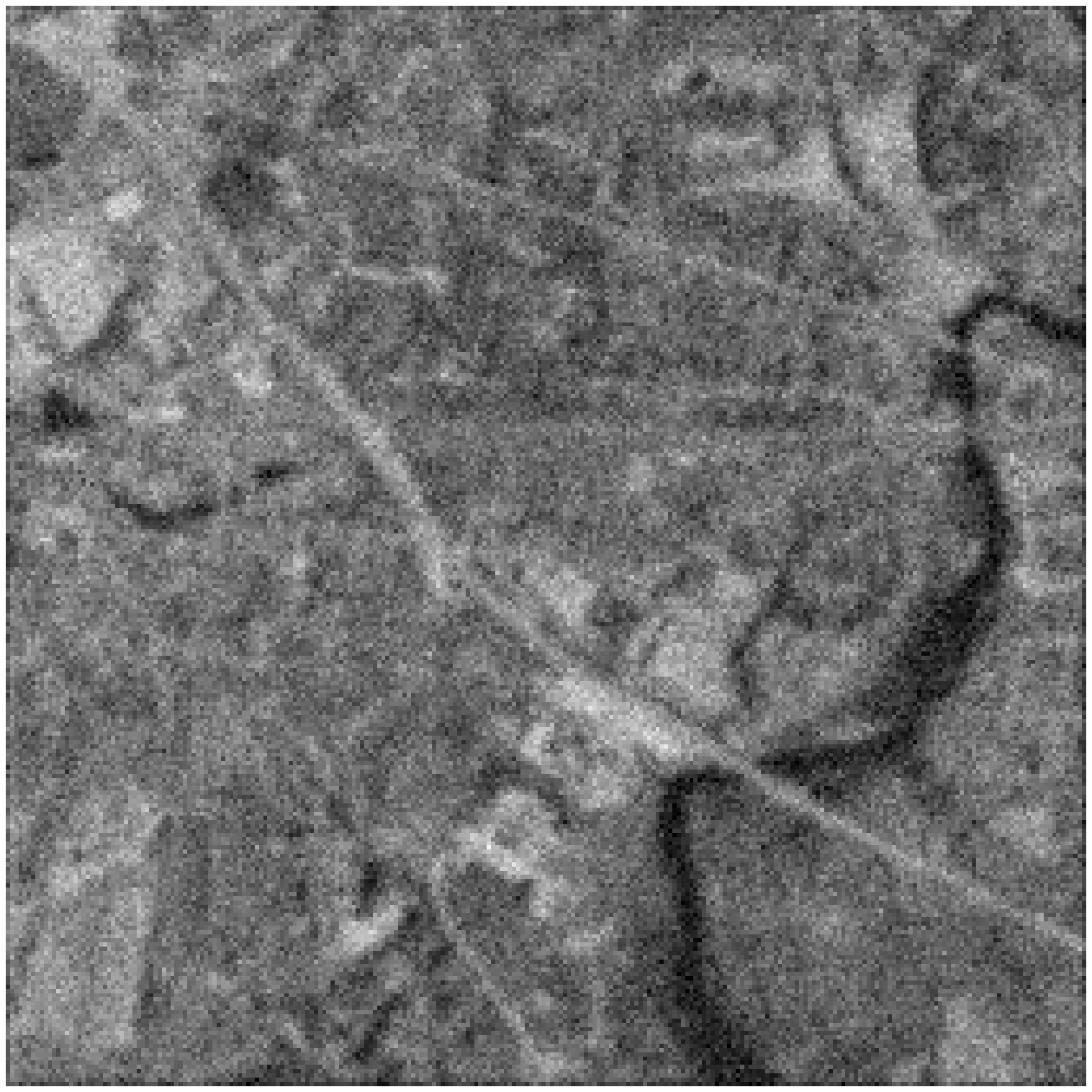} &
 		\includegraphics[scale=0.2]{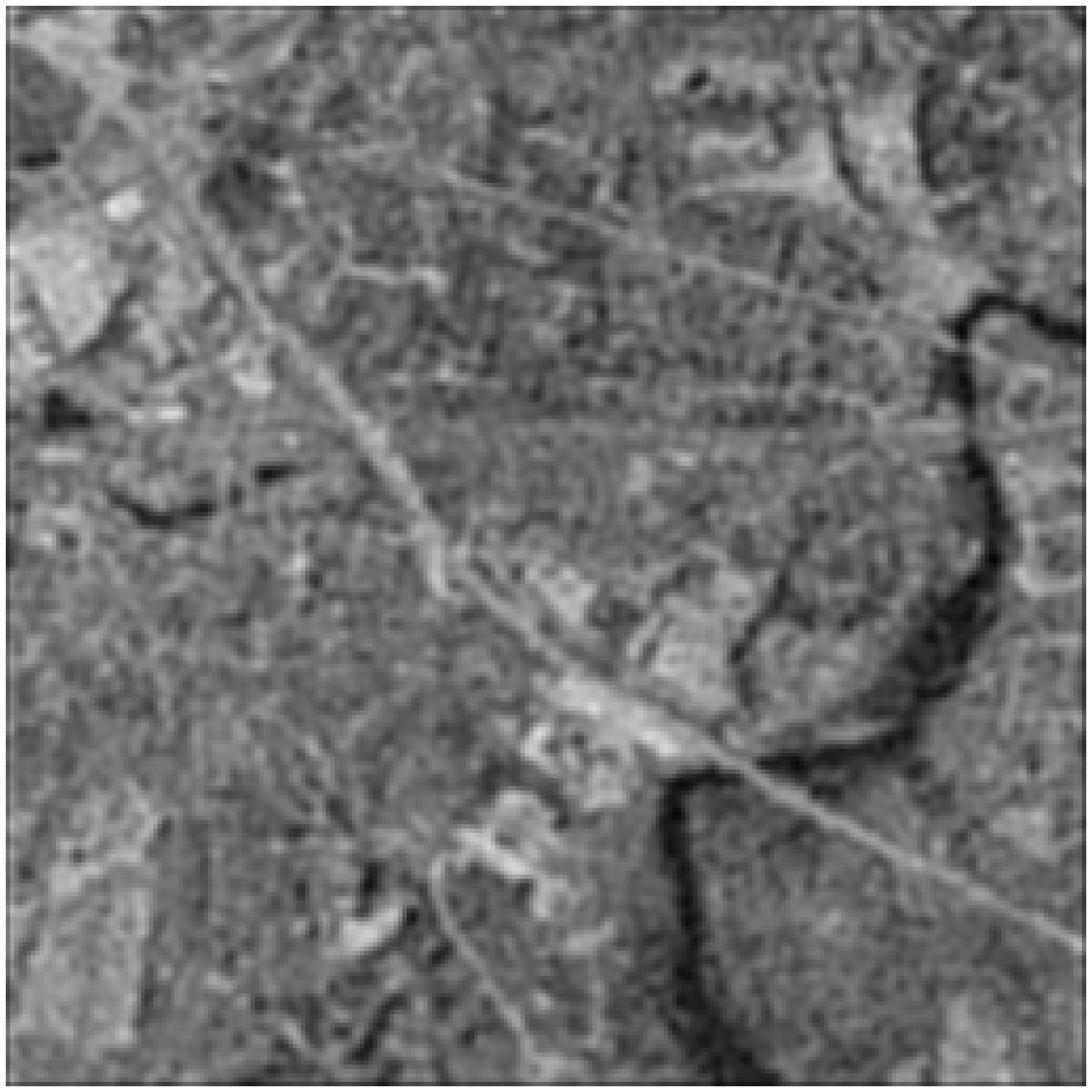} \\
\includegraphics[scale=0.2]{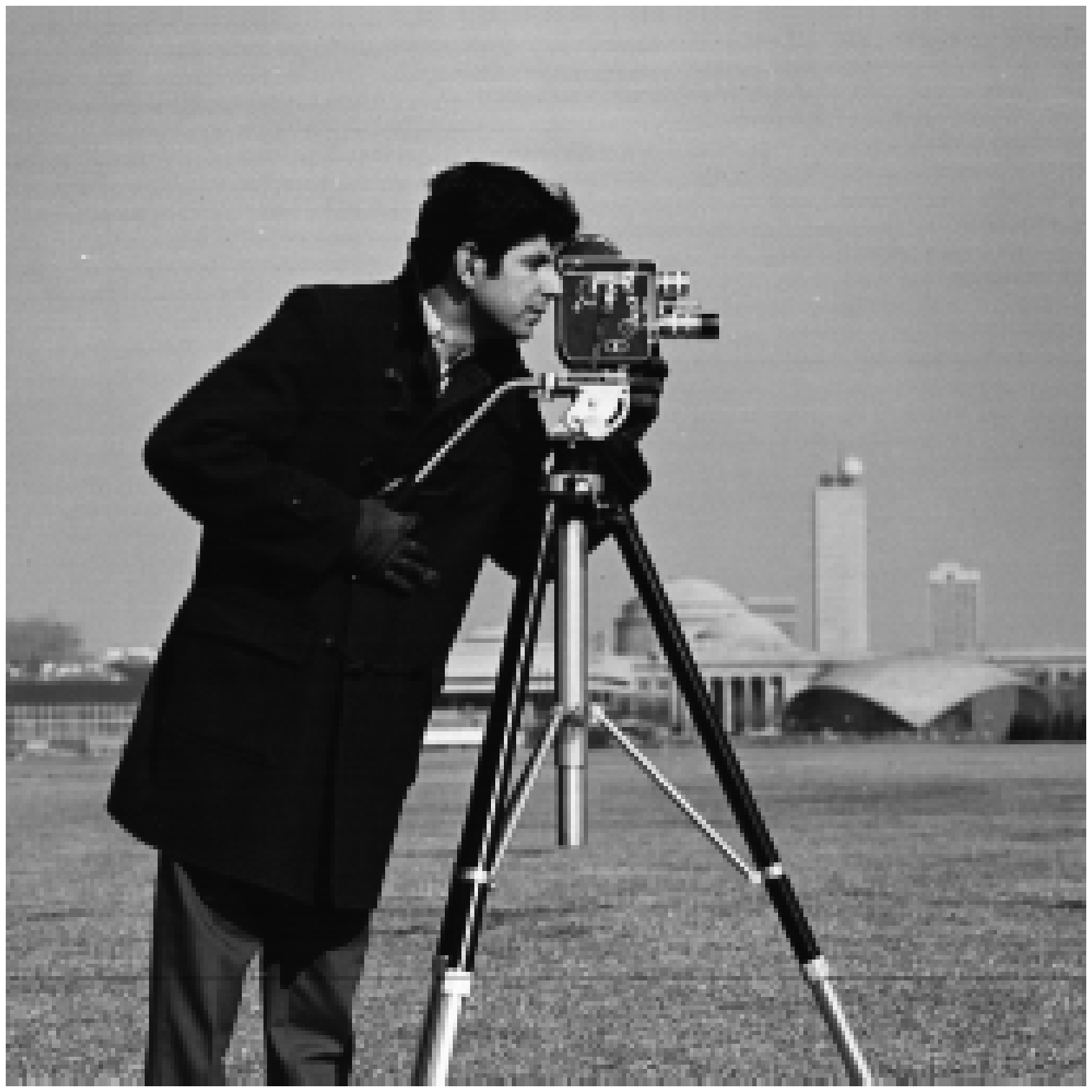} &
 		\includegraphics[scale=0.2]{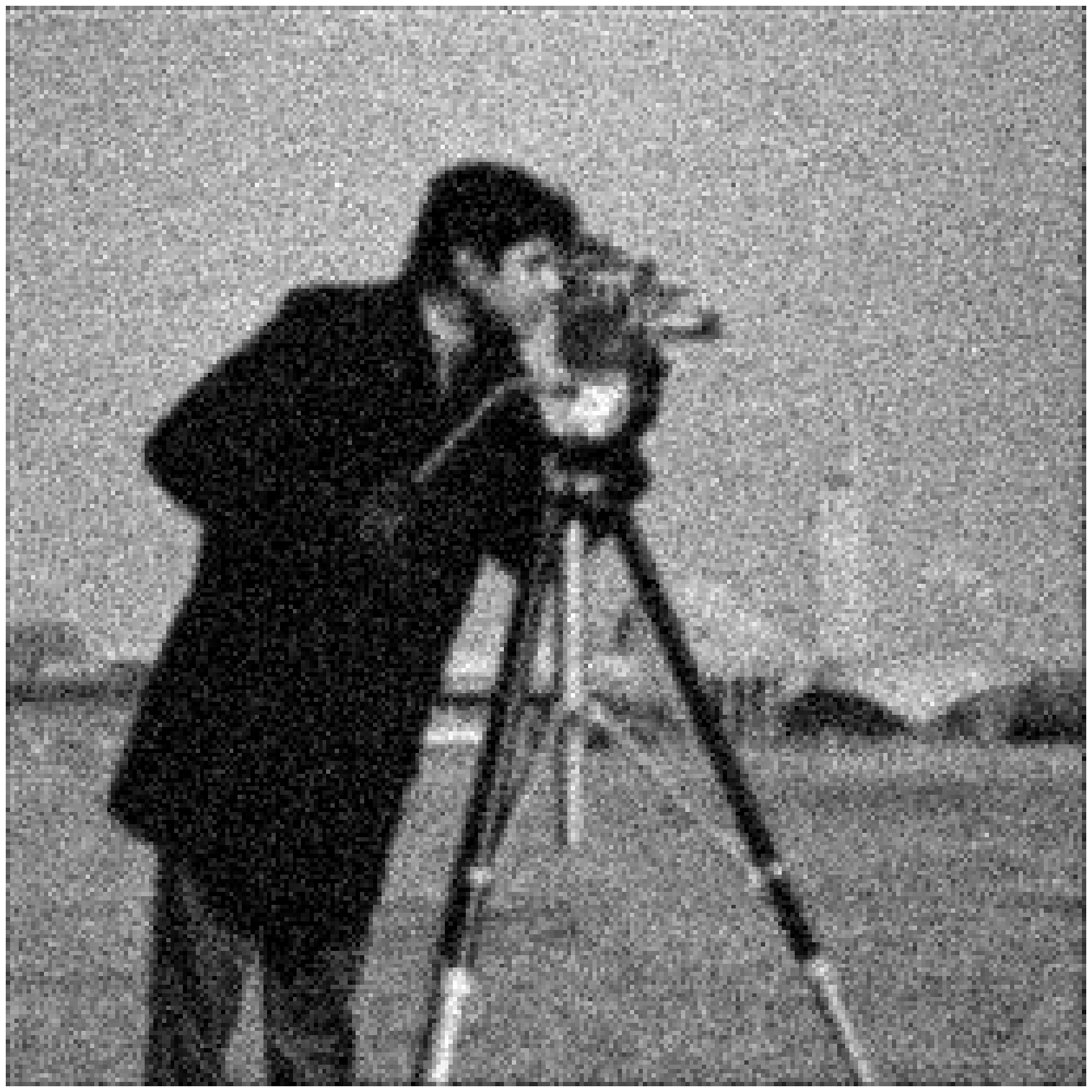} &
 		\includegraphics[scale=0.2]{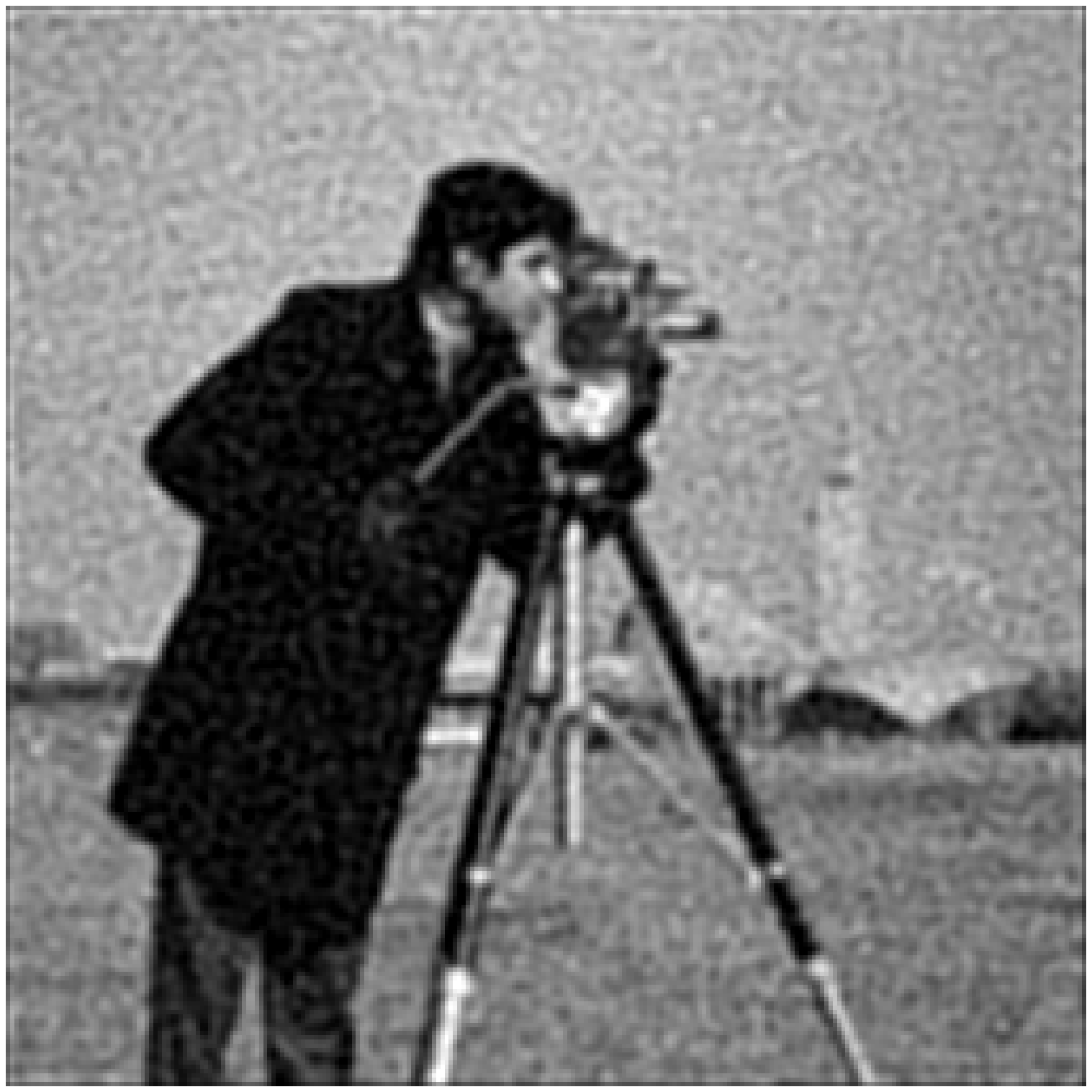}
		\end{tabular}
 	\end{center}		
 \end{figure}

Our empirical results suggest that this is possible. We selected four images from Matlab's {\tt imdemos} image library.  These images were blurred and noise was added.  Then reconstructed images were computed using multi-parameter Tikhonov regularization, using the optimal parameters computed for the training set. Relative reconstruction errors are provided in Table~\ref{tab:MultiLOutside}, and image reconstructions are provided in Figure~\ref{fig:example3_Outside}. To provide a quantitative measure of the similarity between the four images and the set of training images, we report in Table~\ref{tab:MultiLOutside} the average of the 2-norm of the differences and the average of the Mean Structural Similarity (MSSIM) indices between the each of the four images and the set of training images.  MSSIM values that are closer to 1 indicate more structural similarity between images \cite{mssim04}. The reconstruction quality is encouraging, with average reconstruction errors for the ``outsider" images being larger than the median reconstruction error for the validation set by a factor of approximately 2 (we compare to the median of the errors for Tik-multi-$\hat{\lambda}_{GCV}$ in Figure~\ref{fig:example3_box}, but the other two choices would yield similar results). Table~\ref{tab:MultiLOutside} also suggests that, while there is no strict relationship between similarity and image quality, a positive correlation between the two can be expected in general.

\begin{table}[bthp]
	\caption{Comparison of average difference, average SSIM values, and relative reconstruction errors for non-satellite images. Smaller average difference values and average SSIM values closer to 1 imply closer similarity between the image and the set of training images. }
	\label{tab:MultiLOutside} 		
\begin{center}
\scalebox{0.9}{
\begin{tabular}{|l||c|c|c|} \hline	    	
   	&  Ave. Difference   & Ave. SSIM Difference & Rel. Reconstruction Errors \\ \hline
Cell Image    & 71.6595 & 2.348e-01 &  1.13e-02 \\ \hline
Fabric Image  & 83.6700 & 0.482e-01 &  2.50e-02 \\ \hline
Aerial Map    & 90.3935  & 0.432e-01    & 1.75e-02  \\  \hline
Camera Man    & 114.4114 & 1.379e-01   & 1.55e-02  \\  \hline
\end{tabular}}
\end{center}
\end{table}

\section{Conclusions}
\label{sec:conclusions}

In this paper, we described a learning approach to compute optimal regularization parameters for general-form and multi-parameter Tikhonov regularization. We formulated the problem as an empirical Bayes risk minimization problem and developed numerical algorithms for computing the optimal or near-optimal regularization parameters.  Various error measures were considered.  Numerical results illustrate that GSVD filtered solutions give better performance than optimal SVD filtered solution, with less training data.  In addition, optimal parameters for multi-parameter Tikhonov can provide results comparable to GCV, with less computational cost for a large set of images, and these parameters can be used to reconstruct different types of images.

\bibliographystyle{plain}

\begin{thebibliography}{10}

\bibitem{BazanBorgesFrancisco2012}
F.~S.~Viloche Baz{\'a}n, L.~S. Borges, and J.~B. Francisco.
\newblock On a generalization of {R}egi\'nska's parameter choice rule and its
  numerical realization in large-scale multi-parameter {T}ikhonov
  regularization.
\newblock {\em Appl. Math. Comput.}, 219(4):2100--2113, 2012.

\bibitem{BelgeKilmerMiller}
M.~Belge, M.~E. Kilmer, and E.~L. Miller.
\newblock Efficient determination of multiple regularization parameters in a
  generalized {L}-curve framework.
\newblock {\em Inverse Problems}, 18(4):1161--1183, 2002.

\bibitem{bjorck1996numerical}
A.~Bj{\"o}rck.
\newblock {\em Numerical methods for least squares problems}.
\newblock SIAM, Philadelphia, 1996.

\bibitem{BrezinskiMulti2003}
C.~Brezinski, M.~Redivo-Zaglia, G.~Rodriguez, and S.~Seatzu.
\newblock Multi-parameter regularization techniques for ill-conditioned linear
  systems.
\newblock {\em Numer. Math.}, 94(2):203--228, 2003.

\bibitem{ChungChung2013}
J.~Chung and M.~Chung.
\newblock An efficient approach for computing optimal low-rank regularized
  inverse matrices.
\newblock Submitted, 2013.

\bibitem{chung2011designing}
J.~Chung, M.~Chung, and D.~P. O'Leary.
\newblock Designing optimal spectral filters for inverse problems.
\newblock {\em SIAM J. Sci. Comput.}, 33(6):3132--3152, 2011.

\bibitem{ChungKilmerOleary}
J.~Chung, M.~E. Kilmer, and D.~P. O'Leary.
\newblock A framework for regularization via operator approximation.
\newblock Submitted, 2013.

\bibitem{colson2007overview}
B.~Colson, P.~Marcotte, and G.~Savard.
\newblock An overview of bilevel optimization.
\newblock {\em Ann. Oper. Res.}, 153:235--256, 2007.

\bibitem{de2013image}
J.~C. {De los Reyes} and C.~Sch{\"o}nlieb.
\newblock Image denoising: Learning the noise model via nonsmooth
  {P}{D}{E}-constrained optimization.
\newblock {\em Inverse Probl. Imaging}, 7(4), 2013.

\bibitem{GazzolaNovati2013}
S.~Gazzola and P.~Novati.
\newblock Multi-parameter {A}rnoldi-{T}ikhonov methods.
\newblock {\em Electron. Trans. Numer. Anal.}, 40:452--475, 2013.

\bibitem{golubGCV}
G.~H. Golub, M.~Heath, and G.~Wahba.
\newblock Generalized cross-validation as a method for choosing a good ridge
  parameter.
\newblock {\em Technometrics}, 21(2):215--223, 1979.

\bibitem{GolubVanLoanBook}
G.~H. Golub and C.~F.~Van Loan.
\newblock {\em Matrix computations}.
\newblock Johns Hopkins Studies in the Mathematical Sciences. Johns Hopkins
  University Press, Baltimore, third edition, 1996.

\bibitem{haber2003learning}
E.~Haber and L.~Tenorio.
\newblock Learning regularization functionals - a supervised training approach.
\newblock {\em Inverse Problems}, 19(3):611, 2003.

\bibitem{HansenBook1998}
P.~C. Hansen.
\newblock {\em Rank-deficient and discrete ill-posed problems}.
\newblock SIAM Monographs on Mathematical Modeling and Computation. SIAM,
  Philadelphia, 1998.

\bibitem{hansen}
P.~C. Hansen.
\newblock Regularization tools version 4.0 for matlab 7.3.
\newblock {\em Numerical Algorithms}, 46(2):189--194, 2007.

\bibitem{HansenBook2010}
P.~C. Hansen.
\newblock {\em Discrete inverse problems}, volume~7 of {\em Fundamentals of
  Algorithms}.
\newblock SIAM, Philadelphia, 2010.

\bibitem{DeblurringBook}
P.~C. Hansen, J.~G. Nagy, and D.~P. O'Leary.
\newblock {\em Deblurring images: Matrices, spectra, and filtering}, volume~3
  of {\em Fundamentals of Algorithms}.
\newblock SIAM, Philadelphia, 2006.

\bibitem{horesh2009sensitivity}
L.~Horesh and E.~Haber.
\newblock Sensitivity computation of the $\ell_1$ minimization problem and its
  application to dictionary design of ill-posed problems.
\newblock {\em Inverse Problems}, 25(9):095009, 2009.

\bibitem{huang2012optimal}
H.~Huang, E.~Haber, and L.~Horesh.
\newblock Optimal estimation of {$\ell_1$}-regularization prior from a
  regularized empirical {B}ayesian risk standpoint.
\newblock {\em Inverse Probl. Imaging}, 6(3):447--464, 2012.

\bibitem{huber1964robust}
P.~J. Huber.
\newblock Robust estimation of a location parameter.
\newblock {\em The Annals of Mathematical Statistics}, 35(1):73--101, 1964.

\bibitem{kunisch2013bilevel}
K.~Kunisch and T.~Pock.
\newblock A bilevel optimization approach for parameter learning in variational
  models.
\newblock {\em SIAM J. Imaging Sci.}, 6(2):938--983, 2013.

\bibitem{liu2012morphology}
J.~Liu, T.~Liu, L.~de~Rochefort, J.~Ledoux, I.~Khalidov, W.~Chen, J.~Tsiouris,
  C.~Wisnieff, P.~Spincemaille, M.~Prince, and Y.~Wang.
\newblock Morphology enabled dipole inversion for quantitative susceptibility
  mapping using structural consistency between the magnitude image and the
  susceptibility map.
\newblock {\em Neuroimage}, 59(3):2560--2568, 2012.

\bibitem{LuPereverzev2011}
S.~Lu and S.~V. Pereverzev.
\newblock Multi-parameter regularization and its numerical realization.
\newblock {\em Numer. Math.}, 118(1):1--31, 2011.

\bibitem{nagy1998restoring}
J.~G. Nagy and D.~P. O'Leary.
\newblock Restoring images degraded by spatially variant blur.
\newblock {\em SIAM Journal on Scientific Computing}, 19(4):1063--1082, 1998.

\bibitem{peyre2011learning}
G.~Peyr{\'e} and J.~Fadili.
\newblock Learning analysis sparsity priors.
\newblock {\em Sampta'11}, 2011.

\bibitem{Wang2012}
Z.~Wang.
\newblock Multi-parameter {T}ikhonov regularization and model function approach
  to the damped {M}orozov principle for choosing regularization parameters.
\newblock {\em J. Comput. Appl. Math.}, 236(7):1815--1832, 2012.

\bibitem{mssim04}
Z~Wang, A.~C. Bovik, H.~R. Sheikh, and E.~P. Simoncelli.
\newblock Image quality assessment: From error measurement to structural
  similarity.
\newblock {\em IEEE Trans. Image Proc.}, 14, 2004.

\end{thebibliography}

\end{document}